\documentclass{article}

\usepackage{graphicx}

\usepackage{amsmath}
\usepackage{amssymb}
\usepackage{amsthm}
\usepackage[cmtip,arrow]{xy}
\usepackage{pb-diagram,pb-xy}

\newtheorem{theorem}{Theorem}
\newtheorem{prop}[theorem]{Proposition}
\newtheorem{lem}[theorem]{Lemma}
\newtheorem{cor}[theorem]{Corollary}
\theoremstyle{definition}
\newtheorem{remark}{Remark}

\numberwithin{equation}{section}

\numberwithin{theorem}{section}

\newcommand{\ZZ}{{\Bbb Z}}
\newcommand{\RR}{{\Bbb R}}
\newcommand{\CC}{{\Bbb C}}

\newcommand{\PP}{{\Bbb P}}

\newcommand{\Hh}{{\frak H}}

\newcommand{\M}{{\cal M}}

\renewcommand{\O}{{\cal O}}
\newcommand{\T}{{\cal T}}
\newcommand{\R}{{\cal R}}

\renewcommand{\a}{{\mathsf a}}
\renewcommand{\b}{{\mathsf b}}
\renewcommand{\c}{{\mathsf c}}
\renewcommand{\d}{{D}}

\newcommand{\Xmin}[2]{|X^{#1}|\begin{matrix}#2\end{matrix}} 
\renewcommand{\ss}{\scriptscriptstyle}
\newcommand{\sprod}[2]{#1\cdot #2}    
\newcommand{\spisom}[1]{(\sprod{\tilde#1}{\tilde#1})} 
\newcommand{\spbase}[1]{#1\cdot #1} 

\newcommand{\perm}{{S}}  
\newcommand{\tp}[1]{\,{}^{\ss t\!}#1} 
\newcommand{\smsum}{\mathop{\textstyle{\sum}}\limits} 
\newcommand{\im}{\mathop{\rm Im}}

\newcommand{\sgn}{\epsilon} 

\newcommand{\deltadiv}{\theta_\Delta}

\newcommand{\one}{\frak 1}
\newcommand{\two}{\frak 2}
\newcommand{\3}{\frak 3}
\newcommand{\n}{\frak n}
\newcommand{\m}{\frak m}

\newcommand{\kuno}{\kappa}
\newcommand{\kenne}{\kappa}
\newcommand{\kdue}{\kappa}
\newcommand{\supp}{\mathop{\rm supp}}

\newcommand{\Sym}{{\rm Sym}}

\begin{document}
\title{Determinantal
Characterization of Canonical Curves and
Combinatorial Theta Identities}

\author{Marco Matone, Roberto Volpato\\[10pt]
              Dipartimento di Fisica ``G. Galilei'' \\and Istituto
 Nazionale di Fisica Nucleare\\ Universit\`a di Padova,
 Via Marzolo, 8 -- 35131 Padova, Italy}



\date{}

\maketitle

\begin{abstract}\noindent
We characterize genus $g$ canonical curves by the vanishing of
combinatorial products of $g+1$ determinants of Brill-Noether matrices. This also implies the characterization of canonical curves in terms
of $(g-2)(g-3)/2$ theta identities. A remarkable mechanism, based on a basis
of $H^0(K_C)$ expressed in terms of Szeg\"o kernels, reduces such identities
to a simple rank condition for matrices whose entries are logarithmic derivatives of theta functions.
Such a basis, together with the Fay trisecant identity, also leads to the solution of the
question of expressing the determinant of Brill-Noether matrices
in terms of theta functions, without using the problematic Klein-Fay section $\sigma$.\end{abstract}

\section{Introduction}

Let $C$ be a smooth closed Riemann surface of genus $g$;
we will consider non-hyperelliptic
surfaces of genus $g\ge 3$ and we will identify such surfaces with
the corresponding non-singular canonical curves in a projective
space. Recall that a marking for $C$ is given by fixing a standard homotopy basis together with a basepoint
$p_0\in C$.
Let $K_C$ be the canonical line bundle on $C$ and set
$H^0(K_C^n):=H^0(C,K_C^n)$, $n\in\ZZ_{>0}$. Petri's Theorem
\cite{petriuno} determines the ideal of
canonical curves of genus $g\geq4$ by means of relations among elements
of $H^0(K^2_C)$ (see also \cite{ArbHarr,ArbSern,StDon,Schreyer}). As emphasized by Mumford, Petri's relations are
fundamental and should have basic applications (pag. 241 of
\cite{mumfordd}).

Here we show that the $(g-2)(g-3)/2$ linearly
independent relations derived by Petri admit an intrinsic
characterization based on combinatorial properties of determinants,
which in turn imply a characterization of canonical curves in terms
of theta identities. Furthermore, we introduce a remarkable mechanism, based on the choice of a basis
of $H^0(K_C)$ expressed in terms of Szeg\"o kernels, reducing the combinatorial identities
to a simple rank condition for matrices whose entries are logarithmic derivatives of theta functions.
Remarkably, such a basis, together with the Fay trisecant identity, also leads to the solution of the
question of expressing the determinant of the Brill-Noether matrix
$\omega_i(p_j)$ in terms of theta functions, without using the problematic section $\sigma$ considered by Klein and Fay (see
Theorem \ref{brillnoethermatrix} and Remark \ref{import}).

\subsection{Basic notation}

For each $n\in\ZZ_{>0}$, set
$N_n:=h^0(K_C^n)=(2n-1)(g-1)+\delta_{n1}$, $M_n:={g+n-1\choose n}$,
$M:=M_2=g(g+1)/2$, $N:=N_2=3g-3$. Set $I_n:=\{1,\ldots,n\}$ and let $\perm_n$ be
the group of permutations of $n$ symbols.
Let $\{\alpha_1,\ldots,\alpha_g,\beta_1,\ldots,\beta_g\}$ be a
symplectic basis of $H_1(C,\ZZ)$ and  $\{\omega_i\}_{i\in
I_g}$ the basis of $H^0(K_C)$ such that $\oint_{\alpha_i}\omega_j=\delta_{ij}$. Denote by
$\tau_{ij}:=\oint_{\beta_i}\omega_j$ the Riemann period matrix,
$i,j\in I_g$. For each $n\in\ZZ_{>0}$, let $\{\phi_i^n\}_{i\in I_{N_n}}$ be an
arbitrary basis of $H^0(K_C^n)$. We will
denote $\det_{ij}\phi_i^n(x_j)$ by $\det\phi_i^n(x_j)$ or, when
different kind of points need to be specified, by
$\det\phi^n(x_1,\ldots,x_{N_n})$.


\subsection{The strategy}

The first step in our analysis is to define a refinement of original Petri's construction (see \cite{ottimo}).
Choose $p_1,\ldots,p_g\in C$ so that $\det \phi^1(p_1,\ldots,p_g)\ne 0$.
The set $\{\sigma_i\}_{i\in I_g}$, where for all $z\in C$ $$\sigma_i(z):= {\det
\phi^1(p_1,\ldots,p_{i-1},z,p_{i+1},\ldots,p_{g})\over
 \det \phi^1(p_1,\ldots,p_{g})}\ ,
 $$
is a basis of $H^0(K_C)$ adapted to the points $p_1,\ldots,p_g$, i.e. such that $\eta_i(p_j)=0$, for $i\neq j$, $i,j\in I_g$. This definition is independent of the choice of the basis
 $\{\phi^1_i\}_{i\in I_{g}}$ and, up to normalization, of the local coordinates on
 $C$.

This definition of a basis adapted to the points provides a natural way to normalize the holomorphic $1$-differentials. This will be useful when we will need to write polynomials and determinants
of such differentials in terms of theta functions, because such formulae will be automatically independent of the marking.

%
%

\noindent A basic consequence of the above construction is that now
Petri's relations can be expressed in terms of
$\{\omega_i\}_{i\in I_g}$, with respect to which it is possible to
write down a basis for $H^0(K_C^2)$. This is a key point, since the products
$\omega_i\omega_j$ are directly related, via the Kodaira-Spencer
map, to $d\tau_{ij}$. Below we will explicitly show why such a basis and the basis for the
for $H^0(K_C^2)$ are quite natural for our investigation.

\noindent Denote by $\sigma\sigma_i$, $i=1,\ldots,{g+1\choose 2}$, the elements
$\sigma_i\sigma_j$, $i,j\in I_g$, the correspondence being fixed by a bijection $I_M\ni
i\mapsto \one_i\two_i\in \Sym^2(I_g)$.
Consider the
${g+1\choose 2}\times (3g-3)$-dimensional matrix
$\sigma\sigma_i(x_j)$, with $x_1,\ldots,x_{3g-3}\in C$ pairwise
distinct points in general position. One can now
recognize a non-singular submatrix of degree $3g-3$ and, as a consequence,  each one of
the ${g+1\choose 2}-(3g-2)=(g-2)(g-3)/2$ conditions among elements
of $H^0(K_C^2)$ can be expressed in determinantal form as the
vanishing of a minor of degree $3g-2$ containing such a submatrix.
Explicitly, for each one of such relations there is a suitable subset of
$\{\sigma_i\sigma_j\}_{1\le i\le j\le g}$, consisting of $3g-2$
elements, we denote by $\{\sigma\sigma_i\}_{i\in I}$, with $I$ a
suitable set of $3g-2$ indices
such that the relation can be expressed as
$$\det_{\substack{i\in I\\ 1\le j\le
3g-2}}\sigma\sigma_i(x_j)=0
\ , $$
for all $x_1,\ldots,x_{3g-2}\in C$.
Each set $\{\sigma\sigma_i\}_{i\in I}$ spans $H^0(K_C^2)$. As
explained in the Proof of Theorem \ref{main}, such relations can be
simplified and written as the vanishing of minors of degree $2g$.

\noindent
The subsequent step is to express such relations in terms of theta identities. As we will see in Section 4,
there are simple formulae expressing determinants of holomorphic differentials in terms of theta functions and prime forms. However, in order for such formulae to be applied, the determinantal relations above need some non-trivial combinatorial manipulation.
%
The
idea is to express determinants, such as $\det\omega\omega_i(x_j)$,
as a sum of products of determinants of the form $\det\omega_i(x_j)$.
This requires some basic identities, formulated in the
two combinatorial Lemmas, which are of their own interest.

\noindent Let us then illustrate the combinatorics, considered in Sec.\ref{combdetssect}, to express
the determinants of ${g+1\choose
2}$-dimensional matrices
$\omega\omega_i(x_k)$, with $x_k$ arbitrary points of $C$, as polynomials of determinants of the $g\times
g$ Brill-Noether matrices $\omega_i(x_j)$. The starting point is the $g=2$ D'Hoker and Phong remarkable formula \cite{DHokerQP}
$$\left|\begin{matrix}\omega_1^2(x_1) &
\omega_1(x_1)\omega_2(x_1) & \omega_2^2(x_1)\\ \omega_1^2(x_2) &
\omega_1(x_2)\omega_2(x_2) & \omega_2^2(x_2)\\ \omega_1^2(x_3) &
\omega_1(x_3)\omega_2(x_3) & \omega_2^2(x_3)\end{matrix}\right|
=\det\omega(x_1,x_2)\det\omega(x_1,x_3)\det\omega(x_2,x_3) \ .
$$
Although the original proof uses the hyperelliptic form of the abelian differentials, this is actually a purely algebraic identity that holds for any functions $f_1,f_2$. Once the left and right hand side are expressed in terms of theta functions, however, this identity becomes highly non-trivial.
This motivates the search for analogous purely algebraic identities at $g>2$. Namely, for $g$ arbitrary functions
$f_1,\ldots,f_g:A\to{\CC}$ on an arbitrary set $A$, we analyze the
identities between polynomials of determinants of $g\times g$
matrices $\det f_i(x_j)$ and determinants of ${g+1\choose 2}$
dimensional matrices with entries $ff_i(x_k)$, where each $ff_i$
denotes a distinct element in the set of functions
$\{f_if_j\}_{1\leq i\leq j\leq g}$. Such identities are much simpler
than one would naively expect, a striking example being the $g=3$
formula
\begin{align}&\left|\begin{matrix} f_1^2(1) &f_1(1)f_2(1)
&f_1(1)f_3(1) &f_2^2(1) &f_2(1)f_3(1) &f_3^2(1) \notag \\ \vdots & \vdots
& \vdots &\vdots &\vdots &\vdots  \notag \\ f_1^2(6) &f_1(6)f_2(6)
&f_1(6)f_3(6) &f_2^2(6) &f_2(6)f_3(6)
&f_3^2(6)\end{matrix}\right|_{\phantom{\begin{matrix}\vdots\end{matrix}}} \\&\qquad\qquad=\det
f(1,2,3)\det f(1,4,5)\det f(2,4,6)\det f(3,5,6)  \notag \\ &\qquad\qquad\quad-\det
f(6,2,3)\det f(6,4,5) \det f(2,4,1)\det f(3,5,1)\ ,  \notag \end{align}
where $f_i(j):=f_i(x_j)$, which holds for all $x_1,\ldots,x_6\in A$.
Finally, and somehow more surprisingly, it turns out that, under
mild assumptions on the $f_i$ (essentially, $A$
must contain a subset of common zeroes), also particular minors of
the matrices $ff_i(x_j)$ can be expressed as polynomials of
$\det f_i(x_j)$. This result is a priori totally
unexpected and, remarkably, can be directly applied to the
basis of $H^0(K^2_C)$ we consider.

Thanks to the above mentioned combinatorial identities, the above
vanishing conditions can be expressed in terms of polynomials
of determinants of the Brill-Noether matrices, the content of Theorem \ref{main}. Since, as we will show,
each one of such determinants can be directly expressed in
terms of Riemann theta functions, such a construction allows us to
associate a theta relation to each one of the quadrics
generating the ideal of the canonical curve in Petri's construction.
This is the content of Theorem \ref{ththetarel} and Corollary \ref{threla},
whose proofs, together with the one of Theorem \ref{main}, are reported
in Sec.\ref{secthetas}. In Sec.\ref{secthetas} we will introduce bases for $H^0(K_C^n)$ expressed in terms
of Szeg\"o kernels. In the case of $H^0(K_C)$, such a basis, denoted by $\{\lambda_i\}_{i\in I_g}$,
plays a key role in reducing the combinatorial relations among theta functions. This is Theorem \ref{szeghimatr}, proved in Sec.\ref{secthetas},
which provides a simple rank condition for matrices whose entries are
logarithmic derivatives of theta functions. Furthermore, expressing $\det\lambda_i(z_j)$
in terms of theta functions using the Fay trisecant identity \cite{jfayy} and then changing basis to $\{\omega_i\}_{i\in I_g}$,
gives the expression of the determinant of Brill-Noether matrices in terms of theta functions, this is the content of Theorem \ref{brillnoethermatrix}.

Let us further illustrate why the bases $\{\sigma_i\}_{i\in I_g}$ and $\{\sigma\sigma_i\}_{i\in I_N}$ are the natural choices for our
investigation. Suppose one considers an arbitrary basis $\{\eta_i\}_{i\in I_g}$ of $H^0(K_C)$. Under a suitable general position
assumption,
$$\phi_1,\ldots,\phi_{3g-3}:=\eta_1^2,\ldots,\eta_g^2,
\eta_1\eta_2,\ldots,\eta_1\eta_g,\eta_2\eta_3,\ldots,\eta_2\eta_g \ ,$$
will
form a basis of $H^0(K_C^2)$. Then to express any $\eta_i\eta_j$ in terms of
this basis one takes an arbitrary generic collection of points
$q_1,\ldots, q_{3g-3}$ and determines the coefficients $c_{ijk}$ in
$\eta_i\eta_j=\sum c_{ijk}\phi_k$, by requiring the values of both sides to be
equal at $q_1,\ldots, q_{3g-3}$. Note that even if this immediately gives $c_{ijk}$ as a
determinantal quantity in $\eta_i\eta_j(q_l)$ and $\phi_k(q_l)$, a first question with such a choice is how to determine a set of
$M-N=(g-1)(g-2)/2$ linearly independent as above or, at least, as derived by Petri. Actually, without a characterization of the zeroes of the holomorphic differentials, as it happens
both in our and Petri's construction, it does not seem possible to deriving them (unless one does a similar subsequent analysis).
Above we anticipated that using the bases $\{\sigma_i\}_{i\in I_g}$ and $\{\sigma\sigma_i\}_{i\in I_N}$ it allows to apply combinatorial techniques
to find such relations. Furthermore, we are interested in finding relations with well-defined transformation properties under changing of the marking. In particular, a main question is how to
express such relations in terms of the basis $\{\omega_i\}_{i\in I_g}$.
We saw that $\{\sigma_i\}_{i\in I_g}$, and therefore $\{\sigma\sigma_i\}_{i\in I_N}$, is invariant under a change of the marking. In the absence of any specification,
the transformation properties of the basis $\{\eta_i\}_{i\in I_g}$, and therefore of the basis $\{\phi_i\}_{i\in I_N}$, is unknown. Many of the results obtained
in the present paper would be more involved to obtain working with such bases. In particular, using the bases $\{\sigma_i\}_{i\in I_g}$ and $\{\sigma\sigma_i\}_{i\in I_N}$
leads to the $M-N$ determinantal relations involving the $\omega_i$'s. This also leads to the
$M-N$ relations satisfied by $d\tau_i$ which are invariant under a change of the marking. Furthermore, many of the relations between theta functions that we derive
follow by such results.


\subsection{Main results}

The first main result concerns a remarkable expression of the linear relations among holomorphic quadratic differentials
in terms of determinants of holomorphic abelian differentials. The proof is purely algebraic and is based on the combinatorial identities derived in Sec. \ref{combdetssect}.
This is the key step for the expression of such relations in terms of theta identities, as shown in Sec. \ref{secthetas}.

\noindent Set $|x_1,\ldots,x_g|:=
\det\omega_i(x_j)$, for arbitrary $x_1,\ldots,x_g\in C$.

\begin{theorem}\label{main}Let $C$ be a canonical curve of genus $g\geq4$. Fix distinct points
$p_3,\ldots,p_g\in C$ in general position (the precise condition they must satisfy is $K(p_3,\ldots,p_g)\neq0$, with $K$
defined in Eq.\eqref{iltrap}). Then, the ideal of $C$ is generated
by the $(g-2)(g-3)/2$ independent relations
\begin{align}\sum_{s\in\perm_{2g}}&\sgn(s)|x_{s_1},\ldots,x_{s_g}|
|x_{s_g},\ldots,x_{s_{2g-1}}|
|x_{s_1},x_{s_{g+1}},x_{s_{2g}},p_{3},\ldots,\check
p_i,\ldots,p_g| \notag \\ &\cdot |x_{s_2},x_{s_{g+2}},x_{s_{2g}},p_{3},\ldots,\check
p_j,\ldots,p_g|\prod_{k=3}^{g-1}
|x_{s_k},x_{s_{k+g}},p_{3},\ldots,p_g|=0\ ,\notag \end{align}
$3\leq
i<j\leq g$, $s_i:=s(i)$, for all $x_k\in C$, $1\le k\le 2g$, unless
$C$ is trigonal or isomorphic to a smooth plane quintic.\end{theorem}

The second main result concerns remarkable identities among theta functions which follow from linear relations among holomorphic quadratic differentials.

\begin{theorem}\label{szeghimatr}Let $C$ be a canonical curve of genus $g\geq4$, and
choose a symplectic basis of $H_1(C,\ZZ)$. Let
$\delta:=\bigl[{}^{\delta'}_{\delta''}\bigr]\in \ZZ_2^{2g}$ be a
non-singular even theta characteristic and $q\in C$ an arbitrary
point and consider the effective divisor $p_1+\ldots+p_g$ on $C$,
uniquely defined by $A(p_1+\ldots+p_g-q-\Delta)=\delta''+\tau\delta'
\mod \ZZ^g+\tau\ZZ^g$.
Then, the $2g\times g(g-1)/2$ matrix with columns
\begin{equation}\label{szeghione}\left(\begin{matrix}\vec{\nabla} \log\theta[\delta](p_i-p_j)\\ A(p_i-p_j)\end{matrix}\right)_{1\le i<j\le g}\ ,\end{equation}
has rank smaller than $2g-2$.
Futhermore, if
$K(p_3,\ldots,p_g)\neq 0$ and $p_1,\ldots,p_g$ are pairwise
distinct, then the rank of such a matrix is exactly $2g-3$ and, in
particular, the $2g-3$ columns with $i\in\{1,2\}$ and $j\in I_g$, $j>i$, are linearly independent. In
this case, an equivalent statement is that the $g\times (g-1)(g-2)/2$ matrix with columns
\begin{equation}\label{szeghitwo}\left(\begin{matrix}\vec{f}(p_1,p_i,p_j) \end{matrix} \right)_{1<i<j\le g}\ ,\end{equation}  has rank less than $g-1$, where
$$\vec f(p_i,p_j,p_k):=\vec\nabla\log
\theta[\delta](p_{i}-p_{j})+\vec\nabla\log
\theta[\delta](p_{j}-p_{k})+\vec\nabla\log
\theta[\delta](p_{k}-p_{i})
\ ,$$ $i,j,k\in I_g$.  For each non-singular odd
theta characteristic
$\nu:=\bigl[{}^{\nu'}_{\nu''}\bigr]\in\ZZ_2^{2g}$, let $\d_{g-1}$ be
the effective divisor of degree $g-1$ such that
$A(\d_{g-1}-\Delta)=\nu''+\tau\nu'\mod \ZZ^g+\tau\ZZ^g$. If $q$ is in the support of
$\d_{g-1}$, then \begin{equation}\label{szeghithree}{\rm rk}
\left(\begin{matrix}\vec f(p_1,p_2,p_3) & \vec f(p_1,p_2,p_4) & \ldots & \vec f(p_1, p_2,p_g)&
\vec\nabla\theta[\nu](0)\end{matrix} \right)<g-1\ .\end{equation} \end{theorem}

It is worth noticing that the known solution of the
Schottky problem involves non-linear differential equations. In particular,
according to Novikov's conjecture, an indecomposable principally
polarized abelian variety is the Jacobian of a genus $g$ curve if
and only if there exist vectors $U\neq0,V,W\in \CC^g$ such that
$u(x,y,t)=2\partial_x^2\log\theta(Ux+Vy+Wt+z_0,Z)$, satisfies the
Kadomtsev-Petviashvili (KP) equation
$3u_{yy}=(4u_t+6uu_x-u_{xxx})_x$. Relevant progresses on such a
conjecture are due, among the others, to Krichever
\cite{krichever}, Dubrovin \cite{dubrovin} and Mulase \cite{mulase}. Its proof is due to
Shiota \cite{shiota}. A basic step in such a proof concerned the
existence of the $\tau$-function as a global holomorphic function in
the $\{t_i\}$, as clarified by Arbarello and De Concini in
\cite{Arbar}, where it was shown that only a subset of the KP hierarchy
is needed. Their identification of such a subset is based on basic
results by Gunning \cite{gunning}
and Welters \cite{welters,Welterfg},
characterizing the Jacobians by trisecants (see also
\cite{ArbDeConc,donagifirst,donagisecond,donagi,gunningtwo,vangeemen,vansquare}).

The Schottky problem is still under active investigation, see for
example \cite{Sebastian,grushevsky,polishchuk,Marini} for further
developments. In particular, Arbarello, Krichever and Marini proved
that the Jacobians can be characterized in terms of only the first
of the auxiliary linear equations of the KP equation
\cite{ArKrMar,KricheverET}. Very recently appeared the papers \cite{grushsalv} and
 \cite{Kricheverc}  that proved the conjectures by Farkas \cite{Farkas} and Welters \cite{Welterfg}, respectively.

\vskip 6pt

\noindent {\bf Acknowledgments.} We thank the anonymous referee for his/her interest and the constructive review. 

\vskip 6pt

\section{Combinatorial Lemmas for determinants}\label{combdetssect}

Here we first describe a useful formalism for symmetric
tensor products of vector spaces, which will be used throughout the
paper where determinants of holomorphic sections of $K_C^2$ will
play a central role. The first application will be in the next
section where a basis of $H^0(K_C^2)$ in terms of two-fold products
of sections of $H^0(K_C)$ will be introduced. After
describing such a formalism, we will consider the
purely combinatorial problem concerning the determinants of a basis
of a two-fold symmetric product of a finite dimensional space of
functions. In particular, the combinatorial Lemmas \ref{thcombi}
and \ref{thcombvi}
 show that such determinants can be expressed as polynomials in determinants of a basis of the original space of
  functions.
    Such results will be applied to determinantal relations we will derive in the next section, in particular to
  Theorem \ref{thlemma}, leading to Theorem \ref{main}, and then, by Proposition \ref{thdettheta}, expresses them in terms of
  Riemann's theta functions. The latter is the content of Theorem \ref{ththetarel} which, in turn, leads to Corollary \ref{threla},
  whose proofs are reported in Sec.\ref{secthetas}.

\subsection{Identities induced by the isomorphism
$\CC^{M_n}\leftrightarrow \Sym^n(\CC^g)$}\label{notablecombina}

For each $n\in\ZZ_{>0}$, set $I_n:=\{1,\ldots,n\}$ and let $\perm_n$
denote the group of permutations of $n$ elements. Let $V$ be a
$g$-dimensional vector space and let
$M_n:={g+n-1\choose n}$
be the dimension of the $n$-fold symmetric tensor product $\Sym^n(V)$. We denote by
$$\Sym^n(V)\ni\eta_1\cdot\eta_2\cdots\eta_n:=\sum_{s\in\perm_n}\eta_{s(1)}\otimes\eta_{s(2)}\otimes\ldots\otimes
\eta_{s(n)}\ ,$$ the symmetric tensor product of an $n$-tuple
$(\eta_1,\ldots,\eta_n)$ of elements of $V$. Let us fix an isomorphism $\CC^{M_2}\rightarrow\Sym^2(\CC^g)$ and,
more generally, an isomorphism $\CC^{M_n}\rightarrow\Sym^n(\CC^g)$,
$n\in\ZZ_{>0}$.

\begin{remark}\label{definizione}
Let $M:=M_2$, and let $A:\CC^M\rightarrow\Sym^2(\CC^g)$ be the
isomorphism given in coordinates by
$A(\tilde e_i):=\sprod{e_{\one_i}}{e_{\two_i}}$, with  $\{e_i\}_{i\in
I_g}$ and $\{\tilde e_i\}_{i\in I_M}$ the canonical bases of $\CC^g$
and $\CC^M$ respectively, and
$$(\one_i,\two_i):=
\begin{cases}(i,i)\ , &1\le i\le g\ ,\\
(1,i-g+1)\ ,&g+1\le i\le 2g-1\ ,\\(2,i-2g+3)\ , &2g\le i\le 3g-3\
,\\ \hfill\vdots \hfill& \hfill\vdots\hfill\\ (g-1,g)\
,&i=g(g+1)/2\ ,\end{cases}
$$ so that $\one_i\two_i$ is the $i$-th
element in the $M$-tuple $(11,22,\ldots,gg,12,\ldots,1g,23,\ldots)$.
\end{remark}

\noindent
In general, one can define an
isomorphism $A:\CC^{M_n}\to \Sym^n(\CC^g)$, with
$A(\tilde e_i):=e_{\one_i}\cdots e_{\n_i}$, by fixing the $n$-tuples
$(\one_i,\ldots,\n_i)$, $i\in I_{M_n}$, in such a way that $\one_i\le
\two_i\le \ldots\le \n_i$.

\noindent For each vector $u:=\tp(u_1,\ldots,u_g)\in\CC^g$ and matrix $B\in
M_g(\CC)$, set $${u_i^{(n)}=\underbrace{u\cdots
u}_{n\hbox{ times}}}{}_i:=\prod_{\m\in\{\one,\ldots,\n\}}u_{\m_i}\ ,$$
and
$$B^{(n)}_{ij}=(\underbrace{B\cdots B}_{n\hbox{
times}})_{ij}:=\sum_{s\in\perm_n}\prod_{\m\in\{\one,\ldots,\n\}}B_{\m_is(\m)_j}\ ,
$$
$i,j\in I_{M_n}$, where the products
$\prod_{\m\in\{\one,\ldots,\n\}}u_{\m_i}$ and
$\prod_{\m\in\{\one,\ldots,\n\}}B_{\m_is(\m)_j}$ are the standard ones
in $\CC$. In particular, let us define
$$\chi_i^{(n)}:=\prod_{k=1}^g\biggl[\biggl(\sum_{\m\in\{\one,\ldots,\n\}}\delta_{k\m_i}\biggr)!\biggr]=(\delta\cdots\delta)_{ii}\ ,$$ $i\in I_{M_n}$,
(the superscript $(n)$ in $\chi_i^{(n)}$ will be omitted when clear from the context)
where $\delta$ denotes the identity matrix, so that, for example,
$\chi^{(1)}_i=1$, $\chi^{(2)}_i=1+\delta_{\one_i\two_i}$,
$\chi_i^{(3)}=(1+\delta_{\one_i\two_i}+\delta_{\two_i\3_i})(1+\delta_{\one_i\3_i})$.
The coefficient $\chi_i$ is strictly related to the multinomial
coefficient and in particular enters in the multinomial expansion in
$\CC[x_1,\ldots,x_g]$
$$\bigl(\sum_{k=1}^gx_k\bigr)^n=\sum_{i=1}^{M_n}n!\chi_i^{-1}x_{\one_i}\cdots
x_{\n_i}\ .
$$
Such a single indexing satisfies identities repeatedly used in the
following.

\vskip 6pt

\begin{lem}\label{idddd}Let $V$ be a vector space and $f$ an arbitrary function $f:I^n_g\to
V$, where $I^n_g:=I_g\times\ldots\times I_g$. Then, the following
identity holds
\begin{equation}\label{identitt}\sum_{i_1,\ldots,i_n=1}^gf(i_1,\ldots,i_n)=\sum_{i=1}^{M_n}\chi_i^{-1}\sum_{s\in\perm_n}f(s(\one)_i,\ldots,s(\n)_i)\
,\end{equation} that, for $f$ completely symmetric, reduces to
\begin{equation}\label{identittddd}\sum_{i_1,\ldots,i_n=1}^gf(i_1,\ldots,i_n)=n!\sum_{i=1}^{M_n}\chi_i^{-1}f(\one_i,\ldots,\n_i)\
.\end{equation}\end{lem}

\vskip 6pt

\noindent{\sl Proof.} Use
$$\sum_{i_1,\ldots,i_n=1}^gf(i_1,\ldots,i_n)=\sum_{i_n\ge\ldots
\ge i_1=1}^g\sum_{s\in\perm_n}{f(i_{s(1)},\ldots,i_{s(n)})\over
\prod_{k=1}^g[(\sum_{m=1}^n\delta_{ki_m})!] }\ .$$ \hfill$\square$

\vskip 6pt

\noindent Note that for each $u\in\CC^g$, $u^{\otimes n}:=
u\otimes\ldots\otimes u$ is an element of
$\Sym^n(\CC^g)\cong\CC^{M_n}$. By \eqref{identitt}, the following identities
are easily verified
$$u^{\otimes n}\cong \sum_{i=1}^{M_n}\chi_i^{-1}u^{(n)}_i  e^{(n)}_i\ ,\quad (Bu)^{\otimes
n}\cong\sum_{i,j=1}^{M_n}\chi_i^{-1}\chi_j^{-1}B^{(n)}_{ij}
u_j^{(n)} e^{(n)}_i\ ,
$$
where $e_i^{(n)}\in\CC^{M_n}\cong\Sym^n(\CC^g)$, $i\in I_{M_n}$.
Furthermore,
$$\sum_{j=1}^{M_n}\chi_j^{-1}B^{(n)}_{ij}C^{(n)}_{jk}=(BC)^{(n)}_{ik}\ ,$$ where $B,C$ are arbitrary $g\times g$ matrices. For any
non-singular $B$ such an identity yields
$$ \sum_{j=1}^{M_n}\chi_j^{-1}\chi_{k}^{-1}B^{(n)}_{ij}(B^{-1})^{(n)}_{jk}
=\delta^{(n)}_{ik}\chi_k^{-1}=\delta_{ik}\ , $$ and then
$$\det\nolimits_{ij}
\bigl(B^{(n)}_{ij}\chi_j^{-1}\bigr)\det\nolimits_{ij}\bigl((B^{-1})^{(n)}_{ij}\chi_j^{-1}\bigr)
=1\ . $$
In the following, we will denote by
$$B_i:=B_{\one_i\two_i}\ ,
$$
$i\in I_ M$ the elements of a symmetric $g\times g$ matrix.
Furthermore, we will denote the minors of
$B^{(n)}=(B\cdots B)$ by
$$|B\cdots B|^{i_1\ldots i_m}_{j_1\ldots j_m}:=\det_{{i\in \{i_1,\ldots,i_m\}\atop j\in
\{j_1,\ldots,j_m\}}}B^{(n)}_{ij}\ ,$$
$i_1,\ldots,i_m,j_1,\ldots,j_m\in I_{M_n}$, with $m\in I_{M_n}$.

\subsection{Combinatorics of determinants}

Fix a surjection $m:I_g\times I_g\to I_M$, satisfying
$m(i,j)=m(j,i)$, $i,j\in I_g$. It corresponds to an
isomorphism $\CC^M\rightarrow\Sym^2(\CC^g)$ with $\tilde e_{m(i,j)}
\mapsto \sprod{e_i}{e_j}$. For example, by using the construction of
subsec.\ref{notablecombina}, it is possible to define $m$ so that
$m(\one_i,\two_i)=m(\two_i,\one_i)=i$, $i\in I_M$; the corresponding
isomorphism is $A$, introduced in Definition
\ref{definizione}.

\noindent For each map $s:I_M\to I_M$ consider the $g$-tuples $d^k(s)$,
$k\in I_{g+1}$, given by
\begin{equation}\label{ild}d^i_j(s)=d^{j+1}_i(s):=s_{m(i,j)}\ ,\end{equation}
$i\le j$, $i,j\in I_g$. Note that if $s$ is a injective, then
each $g$-tuple consists of distinct integers, and each $i\in I_M$
belongs to two distinct $g$-tuples.

\noindent Consider
$\perm_g^{g+1}:=\underbrace{\perm_g\times\cdots\times\perm_g}_{g+1\hbox{
times}}$ and define $\varkappa:\perm_g^{g+1}\times I_M\rightarrow
I_M$, depending on $m$, by
\begin{equation}\label{ilkappa}\varkappa_{m(i,j)}(r^1,\ldots,r^{g+1})=m(r^i_j,r^{j+1}_i)\
,\end{equation} $i\le j$, $i,j\in I_g$, where
$(r^1,\ldots,r^{g+1})\in\perm_g^{g+1}$. Note that
$$d^i_j(\varkappa(r^1,\ldots,r^{g+1}))=
d^{j+1}_i(\varkappa(r^1,\ldots,r^{g+1}))=m(r^i_j,r^{j+1}_i)\ ,$$
$i\le j$, $i,j\in I_g$. Consider the subset of $I_M$ determined by
$I_{M,n}:=\{m(i,j)|i\in I_n,\,j\in I_g\}$, $n\in I_g$, with the
ordering inherited from $I_M$, and denote by $L:=M-(g-n)(g-n+1)/2$,
its cardinality. The elements $\varkappa_l(r^1,\ldots,r^{g+1})$,
$l\in I_{M,n}$, are independent of $r^j_i$, with $n+1\le i,j\le g$,
and $\varkappa$ can be generalized to a function
$\varkappa:I_{M,n}\times\tilde\perm^{g,n}\to I_M\ ,$ where
$\tilde\perm^{g,n}:=\perm_g^n\times\perm_n^{g-n+1}$, by
\begin{equation}\label{ilkappaii}\varkappa_i(\tilde r^1,\ldots,\tilde r^{g+1}):=
\varkappa_i(r^1,\ldots,r^{g+1})\ ,\end{equation} $i\in I_{M,n}$, $(\tilde
r^1,\ldots,\tilde r^{g+1})\in\tilde\perm^{g,n}$, where
$r^j\in\perm_g$, $j\in I_{g+1}$, are permutations satisfying
$r^j=\tilde r^j$, $j\in I_n$, and $r^j_i=\tilde r^j_i$, $i\in I_n$,
$n+1\le j\le g$. Furthermore, if $\{\varkappa_i(\tilde
r^1,\ldots,\tilde r^{g+1})\}_{i\in I_{M,n}}$ consists of distinct
elements, then it is a permutation of $I_{M,n}$. By a suitable
choice of the surjection \begin{equation}\label{lam}m(j,i)=m(i,j):=M-(g-j)(g-j-1)/2+i\
,\end{equation} $j\le i\in I_g$, we obtain $I_{M,n}=I_L$ as an equality between
ordered sets. Note that this choice for $m$, which is convenient to
keep the notation uncluttered, does not correspond to the
isomorphism introduced in subsec.\ref{notablecombina}.

Consider the maps $s:I\rightarrow I$, where $I$ is any ordered
subset of $I_M$; if $s$ is bijective, then it is a permutation of
$I$. We define the function $\sgn(s)$ to be the sign of the
permutation if $s$ is bijective, and zero otherwise. Let $F$ be a
commutative field and $S$ a non-empty set. Fix a set $\{f_i\}_{i\in
I_g}$, of $F$-valued functions on $S$, and $M$ (possibly coincident)
elements $x_i\in S$, $i\in I_M$. Set $ff_{m(i,j)}:=f_if_j$, $i,j\in
I_g$, and
$$\det f(x_{d^j(s)}):=\det\nolimits_{ik} f_k(x_{d^j_i(s)})\ ,$$
$j\in I_{g+1}$, where $x_i\in S$, $i\in I_M$. Furthermore, for any
ordered set $I\subseteq I_M$, we denote by
$$\det\nolimits_I ff(x_1,\ldots,x_{{\rm Card}(I)})\ ,$$
the determinant of the matrix $(ff_m(x_i))_{^{i\in I_{{\rm
Card}(I)}}_{m\in I}}$.

\begin{lem}\label{thcombi}Choose $n\in I_g$ and $L$
elements $x_i$ in $S$, $i\in I_L$. Fix $g-n$ elements $p_i\in S$,
$n+1\le i\le g$ and $g$ $F$-valued functions $f_i$ on $S$, $i\in
I_g$. The following $g(g-n)$ conditions \begin{equation}\label{combcond}
f_i(p_j)=\delta_{ij}\ ,\end{equation} $1\le i\le j$, $n+1\le j\le g$, imply
\begin{multline}\label{combi}\det\nolimits_{I_{M,n}} ff(x_1,\ldots,x_L)={1\over
c_{g,n}}\sum_{s\in\perm_L}\sgn(s)\prod_{j=1}^n \det
f(x_{d^j(s)})\\ \prod_{k=n+1}^{g+1} \det
f(x_{d^k_1(s)},\ldots,x_{d^k_n(s)},p_{n+1},\ldots,p_g)\end{multline} where
\begin{equation}\label{lacost}c_{g,n}:=\sum_{(\tilde r^1,\ldots,\tilde
r^{g+1})\in\tilde\perm^{g,n}} \prod_{k=1}^{g+1}\sgn(\tilde r^k)
\sgn(\varkappa(\tilde r^1,\ldots,\tilde r^{g+1}))\ .\end{equation} In particular,
for $n=g$ \begin{equation}\label{combii}c_g\det
ff(x_1,\ldots,x_M)=\sum_{s\in\perm_M}\sgn(s) \prod_{j=1}^{g+1} \det
f(x_{d^j(s)})\ ,\end{equation} where
$$
c_g:=c_{g,g}=\sum_{r^1,\ldots,
r^{g+1}\in\perm_g}\prod_{k=1}^g\sgn(r^k)\,\sgn(\varkappa(r^1,\ldots,r^g))
\ .
$$\end{lem}

\vskip 6pt

\noindent {\sl Proof.} It is convenient to fix the surjection $m$ as
in \eqref{lam}, so that $I_{M,n}=I_L$. Next consider \begin{equation}\label{combiii}c_{g,n}
\det\nolimits_{I_L} ff(x_1,\ldots,x_L)=c_{g,n}\sum_{s\in\perm_L}
\sgn(s)ff_1(x_{s_1})\cdots ff_L(x_{s_L})\ ,\end{equation} where $s_i:=s(i)$.
Restrict the sums in \eqref{lacost} to the permutations $(\tilde
r^1,\ldots,\tilde r^{g+1})\in\perm^{g,n}$, $i\in I_n$, such that
$\sgn(\varkappa(\tilde r^1,\ldots,\tilde r^{g+1}))\neq 0$, and set
$s':=s\circ \varkappa(\tilde r^1,\ldots,\tilde r^{g+1})$, so that
$$ff_1(x_{s_1})\cdots
ff_L(x_{s_L})=ff_{\varkappa_1}(x_{s'_1})\cdots
ff_{\varkappa_L}(x_{s'_L})\ ,$$ where $\varkappa_i$ is to be
understood as $\varkappa_i(\tilde r^1,\ldots,\tilde r^{g+1})$. Note
that, for all $l\in I_M$, there is a unique pair $i,j\in I_g$, $i\le
j$, such that $l=m(i,j)$, and by \eqref{ild} and \eqref{ilkappa} the following
identity
$$ff_{\varkappa_l(r^1,\ldots,r^{g+1})}(x_{s'_l})=ff_{m(r^i_j,r^{j+1}_i)}(x_{s'_{m(i,j)}})=
f_{r^i_j}(x_{d^i_j(s')}) f_{r^{j+1}_i}(x_{d^{j+1}_i(s')})\ ,$$ holds
for all $(r^1,\ldots,r^{g+1})\in\perm_g^{g+1}$. On the other hand,
if $l\in I_L$, then $i\le n$ and by \eqref{ilkappaii}
\begin{align}f_1(x_{s_1})\cdots ff_L(x_{s_L})
&=\prod_{i=1}^n f_{\tilde r^i_1}(x_{d^i_1(s')})\cdots f_{\tilde
r^i_g} (x_{d^i_g(s')}) \notag \\
&\cdot\prod_{j=n+1}^{g+1} f_{\tilde
r^j_1}(x_{d^j_1(s')})\cdots f_{\tilde r^j_n} (x_{d^j_n(s')})\ .  \notag \end{align}
The condition $f_i(p_j)=\delta_{ij}$, $i\le j$, implies
$$\sum_{\tilde r^j\in\perm_n}\sgn(\tilde r^j)f_{\tilde r^j_1}(x_{d^j_1(s')})
\cdots f_{\tilde r^j_n}(x_{d^j_n(s')})= \det
f(x_{d^j_1(s')},\ldots,x_{d^j_n(s')},p_{n+1},\ldots,p_g)\ ,$$
$n+1\le j\le g+1$. Hence, \eqref{combi} follows by replacing the sum over
$s$ with the sum over $s'$ in \eqref{combiii}, and using
$\sgn(s)=\sgn(s')\,\sgn (\varkappa(\tilde r^1,\ldots,\tilde
r^{g+1}))$. Eq.\eqref{combii} is an immediate consequence of
\eqref{combi}.\hfill$\square$

\vskip 6pt

\begin{remark}\label{rrremark}The summation over $\perm_M$ in \eqref{combii}
yields a sum over $(g+1)!$ identical terms, corresponding to
permutations of the $g+1$ determinants in the product. This
overcounting can be avoided by summing over the following subset of
$\perm_M$
$$\perm'_M:=\{s\in \perm_M, \hbox{ {\it s.t.} } s_1=1,\ s_2<s_3<\ldots<s_g,\
s_2<s_i,\ g+1\leq i\leq 2g-1\}\ ,$$ and by replacing $c_g$ by
$c_g/(g+1)!$.\end{remark}

\vskip 6pt

\noindent
Direct computation gives \begin{equation}\label{glic}c_{g,1}=g!\ , \;
c_{g,2}=g!(g-1)!(2g-1)\ , \; c_2=6\ , \; c_3=360\ , \;
c_4=302400\ .\end{equation} For $g=2$, $c_g/(g+1)!=1$ and
$\perm'_{M=3}=\{(1,2,3)\}$, so that \begin{equation}\label{peppp} \det
ff(x_1,x_2,x_3)=\det f(x_1,x_2)\det f(x_1,x_3)\det f(x_2,x_3)\ .\end{equation}

\noindent Remarkably, Eq.\eqref{combii} simplifies considerably for $g=3$. In this
case, $M=6$ and the identity considered in the Introduction holds.
Such a formula can be proved by direct computation; it can also be
derived, up to a constant, by noting that the right hand side is
completely antisymmetric in $x_1,\ldots,x_6$.

A crucial point in proving Lemma \ref{thcombi} is that if
$\varkappa_i(\tilde r^1,\ldots,\tilde r^{g+1})$, $i\in I_{M,n}$, are
pairwise distinct elements in $I_M$, then they belong to
$I_{M,n}\subseteq I_M$, with $\varkappa$ a permutation of such an
ordered set. For a generic ordered set $I\subseteq I_M$, one should
consider $\varkappa$ as a function over $g+1$ permutations $\tilde
r^i$, $i\in I_{g+1}$, of suitable ordered subsets of $I_g$. In
particular, $\tilde r^i$ should be a permutation over all the
elements $j\in I_g$ such that $m(i,j)\in I$, for $j\ge i$, or
$m(i-1,j)\in I$, for $j<i$. However, the condition that the elements
$\varkappa_i(\tilde r^1,\ldots,\tilde r^{g+1})$, $i\in I$, are
pairwise distinct does not imply, in general, that they belong to
$I$ and Lemma \ref{thcombi} cannot be generalized to a determinant of
products $ff_i$, $i\in I$. On the other hand, the subsets
\begin{equation}\label{laI}I:=I_{M,n}\cup \{m(i,j)\}\ ,\end{equation} satisfy such a condition for
$n< i,j\le g$ and yield the following generalization of Lemma
\ref{thcombi}.

\vskip 6pt

\begin{lem}\label{thcombvi}Assume the hypotheses of Lemma {\rm \ref{thcombi}} for $n<g$,
and choose an element $x_{L+1}\in S$, and a pair $i,j$, $n< i,j\le
g$. Then the following relation
\begin{align}\label{combvi}\det\nolimits_I
ff(x_1,\ldots,x_{L+1})&={1\over
c'_{g,n}}\sum_{s\in\perm_{L+1}}\sgn(s) \prod_{k=1}^n\det
 f(x_{d^k(s)}) \notag \\ &\cdot\det
f(x_{d^{n+1}_1(s)},\ldots,x_{d^{n+1}_{n+1}(s)},p_{n+1},\ldots,\check
p_i,\ldots,p_g)\notag\\ &\cdot \det
f(x_{d^{n+2}_1(s)},\ldots,x_{d^{n+2}_{n+1}(s)},p_{n+1},\ldots,\check
p_j,\ldots,p_g) \notag\\&\cdot\prod_{l=n+3}^{g+1}\det
f(x_{d^l_1(s)},\ldots,x_{d^l_n(s)},p_{n+1},\ldots,p_g) \ ,\notag\end{align}
where
$$c'_{g,n}:=\sum_{(\tilde r^1,\ldots,\tilde r^{g+1})\in\tilde\perm^I}\prod_{i=1}^{g+1}\sgn(\tilde
r^i)\sgn(\varkappa(\tilde r^1,\ldots,\tilde r^{g+1}))\ ,$$
$\tilde\perm^I:=\perm_g^n\times\perm_{n+1}^2\times\perm_n^{g-n-1}$,
with $I$ defined in {\rm \eqref{laI}}, holds.\end{lem}

\vskip 6pt

\noindent {\sl Proof.} A straightforward generalization of the proof
of Lemma \ref{thcombi}. \hfill$\square$

\section{Relations in $H^0(K_C^2)$}\label{primecostr}

In this section we introduce a distinguished basis of $H^0(K_C)$ on
a compact non-hyperelliptic Riemann surface $C$ of genus $g$, which
we will identify with its corresponding canonical curve in
$\PP^{g-1}$. We consider a suitable set of products of pairs of such
abelian differentials that, under some conditions we will discuss
in detail, corresponds to a basis of $H^0(K_C^2)$. This construction
allows us to explicitly derive $(g-2)(g-3)/2$ linear relations among
elements of $H^0(K_C^2)$ corresponding to a set of generators for
the degree $2$ ideal of the canonical curve.

The main novelties in the present approach are the
normalization of Petri's basis and the expression of
quadrics generating the ideal of the canonical curve in
terms of determinantal relations given in Theorem \ref{thlemma}.
These results will be used to derive
some of the results of Sec.\ref{secthetas}, such as Theorem \ref{main}, which follows
by applying the combinatoric Lemma \ref{thcombvi} to Theorem \ref{thlemma}.

\subsection{A distinguished basis of $H^0(K_C^2)$}

Let $C$ be a canonical curve of genus $g$. For $n\in\ZZ$, denote by $C_n=\Sym^n(C)$ the space of effective divisors of degree $n$ on $C$.
Set $\eta_i:=\phi^1_i$, $i\in I_g$, and fix the divisor $\c:=x_1+\ldots+x_{g-1}$ in such a
way that the matrix $[\eta_i(x_j)]_{^{i\in I_g}_{j\in I_{g-1}}}$ be
of maximal rank. $\sigma_{c}(p,q):={\det
\eta(p,x_1,\ldots,x_{g-1})/
 \det \eta(q,x_1,\ldots,x_{g-1})}$, is a meromorphic function on
$C_{g-1}$ and a meromorphic section of the bundle ${\cal
L}:=\pi^*_1K_C\otimes \pi^*_2 K_C^{-1}$ on $C\times C$ that is identically one on the diagonal. Here,
$\pi_1$, $\pi_2$ denote the projections of $C\times C$ onto its first
and second component, respectively.

\vskip 6pt

\begin{prop}\label{thnewbasis}Fix $n\in\ZZ_{>0}$ and let $p_1,\ldots,p_{N_n}$ be a set of points
of $C$ such that $\det \phi^n(p_1,\ldots,p_{N_n})\ne 0$, with
$\{\phi^n_i\}_{i\in I_{N_n}}$ an arbitrary basis of $H^0(K_C^n)$.
Then \begin{equation}\label{basendiff}\gamma^n_i(z):= {\det
\phi^n(p_1,\ldots,p_{i-1},z,p_{i+1},\ldots,p_{N_n})\over
\det \phi^n(p_1,\ldots,p_{N_n})}\ ,\end{equation} $i\in I_{N_n}$, for all $z\in C$, is a basis of
$H^0(K_C^n)$ which is independent of the choice of the basis
$\{\phi^n_i\}_{i\in I_{N_n}}$ and, up to normalization, of the local coordinates on
$C$.\end{prop}

\vskip 6pt

\noindent This fixes the bases of
$H^0(K_C^n)$ adapted to a $N_n$-tuple of points $p_1,\ldots,p_{N_n}$, i.e. 
$\gamma_i^n(p_j)=\delta_{ij}$, $i,j\in I_{N_n}$, and
$\det\gamma_i^n(p_1,\ldots,p_{j-1},z,p_{j+1},\ldots,p_{N_n})=\gamma_j^n(z)$
for all $z\in C$.
Furthermore, $\det \gamma_i^n(z_j)={ \det
\phi_i^n(z_j)\over\det \phi_i^n(p_j)}$ for all $z_1,\ldots,z_{N_n}\in
C$.
In the case $n=1$ the basis adapted to $p_1,\ldots,p_g$ will appear frequently in our investigation. In particular, a choice of $p_1,\ldots,p_g\in C$ such that $\det
\eta_i(p_j)\ne 0$ determines a basis of $H^0(K_C)$, given by, $i\in
I_g$, \begin{equation}\label{newbasis}\sigma_i(z):=\gamma^1_i(z)\ .\end{equation}

\noindent Each basis $\{\eta_i\}_{i\in I_g}$ of $H^0(K_C)$ naturally
defines a basis in $\Sym^2(H^0(K_C))$ by
$\sprod{\eta}{\eta}_i:=\sprod{\eta_{\one_i}}{\eta_{\two_i}}$, $i\in
I_M$. Set
\begin{equation}\label{lev}v_i:=\psi(\spbase{\sigma}_i)=\sigma_{\one_i}\sigma_{\two_i}\ ,\end{equation}
for $i\in I_M$, where $\psi$ is the natural map
$\psi\colon\Sym^2(H^0(K_C))\to H^0(K_C^2)$ and note that
\begin{equation}\label{devv}v_i(p_j)=
\begin{cases}\delta_{ij}\ , & i\in I_g \
,\\ 0\ ,& g+1\le i\le M\ ,
\end{cases}
\end{equation} $j\in I_g$. By Max
Noether's Theorem, if $C$ is a Riemann surface of genus two or
non-hyperelliptic with $g\ge 3$, then $\Sym^2(H^0(K_C)) \twoheadrightarrow H^0(K_C^2)$. By
dimensional reasons, it follows that for $g=2$ and $g=3$ in the
non-hyperelliptic case, the set $\{v_i\}_{i\in I_N}$ is a basis of
$H^0(K_C^2)$ if and only if the basis adapted to $p_1,\ldots,p_g$, that is $\{\sigma_i\}_{i\in I_g}$, is a basis of
$H^0(K_C)$. On the other hand, for $g\ge 3$ in the hyperelliptic
case the map $\psi$ is not surjective, so
that $v_1,\ldots,v_N$ are not linearly independent. The other
possibilities are considered in the following proposition whose proof is standard
(see, for example, \cite{ottimo,StDon}).

\vskip 6pt

\begin{prop}\label{thlev}Fix the points $p_1,\ldots,p_g\in C$, with $C$ non-hyperelliptic of genus
$g\ge 4$. If the following conditions are satisfied
\smallskip
\item{\it i.} $\det \eta_i(p_j)\ne 0$, with $\{\eta_i\}_{i\in I_g}$ an arbitrary basis of $H^0(K_C)$;
\item{\it ii.} the divisors $(\sigma_1)-\b$ and $(\sigma_2)-\b$, where $\b:=\sum_{i=3}^gp_i$ and $\{\sigma_i\}_{i\in I_g}$ is
defined in \eqref{newbasis}, have disjoint supports;\\
\smallskip
\noindent then $\{v_i\}_{i\in I_N}$ is a basis of $H^0(K_C^2)$.\end{prop}

\subsection{Relations in $H^0(K_C^2)$}

By Max Noether's theorem, the
natural map $\psi:\Sym^2(H^0(K_C))\rightarrow H^0(K_C^2)$ is surjective if $C$ is canonical. In this section, an explicit description of the kernel of $\psi$ will be given, depending on the choice of a basis of $H^0(K_C)$.

\noindent Let us first fix our notation. Given a basis  $\{\phi^n_i\}_{i\in I_{N_n}}$ of $H^0(K_C^n)$, let us denote by $\tilde\phi^n:H^0(K_C^n)\rightarrow \CC^{N_n}$ the
isomorphism such that $\tilde\phi^n(\phi^n_i)=e_i$, with $\{e_i\}_{i\in
I_{N_n}}$ the canonical basis of $\CC^{N_n}$, and by $\spisom{\phi^n}:\Sym^2(H^0(K_C^n))\rightarrow\Sym^2(\CC^{N_n})$ the induced isomorphism of symmetric product spaces.
Consider $\{\sigma_i\}_{i\in I_g}$, that is the basis of $H^0(K_C)$ adapted to $p_1,\ldots,p_g$ we constructed in
the previous subsection. In the following, we will derive explicitly the
matrix form of the map $\tilde v\circ\psi\circ\spisom{\sigma}^{-1}$
$$
\CC^M\cong \Sym^2(\CC^g)\xrightarrow{\makebox[35pt][c]{{\small $\spisom{\sigma}^{-1}$}}}\Sym^2(H^0(K_C))\xrightarrow{\makebox[35pt][c]{{\small $\psi$}}}H^0(K_C^2)\xrightarrow{\makebox[35pt][c]{{\small $\tilde v$}}}\CC^M
$$
and this will lead to the explicit expression
of $\ker\psi$.

For each set $\{\phi_i^n\}_{i\in I_{N_n}}\subset H^0(K_C^n)$,
consider
the Wronskian $W[\phi^n](p):=\det \partial_p^{j-1}\phi^n_i(p)$. If, for some $p_1,\ldots,p_{N_n}\in C$,
$\det \phi^n_i(p_j)$ does not vanish identically, then,
for each $\{\phi_i^{n'}\}_{i\in I_{N_n}}\subset H^0(K_C^n)$, we have the constant
$$\R_{{\phi}^{n}}[\phi^{n'}]:={\det \phi^{n'}_i(p_j)\over
\det \phi^n_i(p_j)}={W[\phi^{n'}](p)\over
W[{\phi}^{n}](p)}\ .
$$

\subsection{The ideal of $C$ as determinantal variety}

We now show that the ideal of a canonical curve $C$ is generated by determinantal conditions. Besides its own interest, such a result, together with the combinatorics developed in the previous section, will
allow to express such determinantal relations in terms of determinants of the Brill-Noether matrices $\omega_i(x_j)$. In this way the determinants of the holomorphic quadratic differentials reduce to expressions involving
only determinants of the holomorphic abelian differentials. This is the content of Theorem \ref{main}. In turn, this also resolves the problem of expressing the determinantal conditions in terms of theta functions. Actually, whereas the determinantal relations for the quadratic differentials, derived in the next theorem, cannot be directly expressed in terms of theta functions, this can be done once such relations are expressed in terms of determinants of the Brill-Noether matrices. In this way
one can then use their expression in terms of theta functions considered in the next section and then express the determinantal relations as theta relations.

\begin{theorem}\label{thlemma}The ideal of a canonical curve $C$ of genus $g\geq4$
is generated
by the $(g-2)(g-3)/2$ determinantal conditions
\begin{equation}\label{laprincipale}\det\nolimits_I \sigma\sigma(x_1,\ldots,x_{2g})=0
\ ,\end{equation} where $$I:=I_{M,2}\cup\{m(i,j)\}=\{m(1,1),\ldots,m(1,g),m(2,2),\ldots,m(2,g),m(i,j)\}\ ,$$ $3\leq i<j\leq g$, for all $x_i\in
C$, $i\in I_{2g}$ unless $C$ is trigonal or isomorphic to a smooth
plane quintic.\end{theorem}

\vskip 6pt

\noindent {\sl Proof.} For all $i\in I_N$, $j\in I_M$ set $
\tilde\psi_{ij}:=\R_{v}[v_1,\ldots,v_{i-1},v_j,v_{i+1},\ldots,v_N]$.  $v_1,\ldots,v_M$ satisfy the following relations
\begin{equation}\label{lemma}v_i=\sum_{j=1}^N\tilde\psi_{ji}v_j=\sum_{j=g+1}^N\tilde\psi_{ji}v_j\
,\end{equation} $i=N+1,\ldots,M$, which follows by
Cramer's rule and by the identities \eqref{devv}
that imply $\tilde\psi_{ji}=0$ for $j\in I_g$ and $i=N+1,\ldots,M$.
Eq.\eqref{laprincipale} is equivalent to \eqref{lemma}. The rest of the theorem
then follows by Petri's Theorem. \hfill$\square$

\vskip 6pt

\noindent In the following we will show how from \eqref{laprincipale} or, equivalently, \eqref{lemma}, one can get the explicit expression of the coefficients of the $M-N$ quadrics
containing a canonical curve. Even if to derive Eq.\eqref{laprincipale} we used the particular bases of $H^0(K_C)$ and $H^0(K^2_C)$, by Eq.\eqref{laprincipale} one may derive
such relations involving an arbitrary basis $\{\eta_i\}_{i\in I_g}$ of $H^0(K_C)$ as well. In section 4 we will express the coefficients of such quadrics in terms of theta functions in the
relevant case when the basis of $H^0(K_C)$ is $\{\omega_i\}_{i\in I_g}$.

\noindent
Let $\iota:\CC^N\rightarrow \CC^M$ be the injection
$\iota(e_i)=\tilde{e}_i$, $i\in I_N$. Consider the map
$\iota\circ\tilde\psi:\CC^M\rightarrow \CC^M$, where
$\tilde\psi:\CC^M\rightarrow \CC^N$ is the homomorphism with matrix
elements $\tilde\psi_{ij}$, so that
$$(\iota\circ\tilde\psi)_{ij}=
\begin{cases}\tilde\psi_{ij}\ , &1\le i\le N\ ,\\ 0\ ,&N+1\le
i\le M\ ,
\end{cases}
$$ $j\in I_M$. Since
$(\iota\circ\tilde\psi)_{ij}=\delta_{ij}$, $\forall\,i,j\in I_N$, we
have
$\sum_{i=1}^M(\iota\circ\tilde\psi)_{ji}(\iota\circ\tilde\psi)_{ik}=
\sum_{i=1}^N(\iota\circ\tilde\psi)_{ji}\tilde\psi_{ik}=(\iota\circ\tilde\psi)_{jk}$,
$j,k\in I_M$. Hence, $\iota\circ\tilde\psi$ is a projection of rank
$N$ and, since $\iota$ is an injection,
$\ker\tilde\psi=\ker \iota\circ\tilde\psi= ({\rm
id}-\iota\circ\tilde\psi)(\CC^M)$.

\noindent Standard arguments imply the following lemma.

\vskip 6pt

\begin{lem}\label{thutilde} The set $\{\tilde u_{N+1},\ldots,\tilde u_M\}$, $\tilde u_i:=\tilde e_i-\sum_{j=1}^N\tilde
e_j\tilde\psi_{ji}$, $N+1\leq i\le M$, is a basis of
$\ker\tilde\psi$.\end{lem}

\vskip 6pt

\noindent Set $\eta\eta_i:=\psi(\spbase{\eta}_i)$, $i\in I_g$, and
let $X^\eta$ be the automorphism on $\CC^M$ in the commutative
diagram
$$\begin{diagram} \node{\Sym^2(H^0(K_C^2))} \arrow{s,l}{ \spisom{\sigma}} \arrow{se,t}{\spisom{\eta}}\\
\node{\CC^M}\arrow{e,t}{X^\eta}\node{\CC^M}
\end{diagram}
$$
%
whose matrix elements are
\begin{equation}\label{XXXX}X^\eta_{ji}=\chi_j^{-1}([\eta]^{-1}[\eta]^{-1})_{ij}={[\eta]^{-1}_{\one_i\one_j}[\eta]^{-1}_{\two_i\two_j}+
[\eta]^{-1}_{\one_i\two_j}[\eta]^{-1}_{\two_i\one_j}\over
1+\delta_{\one_j\two_j}}\ ,\end{equation} $i,j\in I_M$, where $[\eta]$ is the matrix
with entries $[\eta]_{ij}:=\eta_i(p_j)$, so that
\begin{equation}\label{vXeta}\sigma\sigma_i=\sum_{j=1}^MX^\eta_{ji}\,\eta\eta_j\ ,\end{equation} $i\in I_M$.
Since $\eta\eta_i$, $i\in I_M$, are linearly dependent, the matrix
$X^\eta_{ij}$ is not uniquely determined by \eqref{vXeta}.

\vskip 6pt

\begin{theorem}\label{thcorol}
\begin{equation}\label{corol}\sum_{j=1}^MC^\eta_{ij}\eta\eta_j=0\ ,\end{equation} $N+1\leq i\le
M$, where \begin{equation}\label{leC} C^\eta_{ij}:=\sum_{k_1,\ldots,k_N=1}^M
\Xmin{\eta}{\ss 1 \hfill \ldots \hfill N i \\ \ss k_1\hfill \ldots
\hfill k_N j\hfill}\; \R_v[\eta\eta_{k_1},\ldots,\eta\eta_{k_N}] \
,\end{equation} are $M-N$ independent linear relations among elements of
$H^0(K_C^2)$. Furthermore, for all $p\in C$
\begin{equation}\label{detv}W[v](p)=\sum_{i_1,\ldots,i_N=1}^M\Xmin{\eta}{\ss 1 \hfill
\ldots \hfill N \\ \ss i_1\hfill \ldots \hfill
i_N}\;W[\eta\eta_{k_1},\ldots,\eta\eta_{k_N}](p)\ .\end{equation}\end{theorem}

\vskip 6pt

\noindent {\sl Proof.} Eq.\eqref{corol} follows by \eqref{lemma}, \eqref{vXeta}
and the identity
$$\sum_{l=1}^N(-)^lX^\eta_{jl}\;\Xmin{\eta}{\ss i\,1\hfill\ldots \check{l}
\ldots \hfill N\\ \ss k_1\hfill\ldots \ldots\hfill k_N}+
X^\eta_{ji}\; \Xmin{\eta}{\ss 1\hfill\ldots\hfill N\\ \ss
k_1\hfill\ldots\hfill k_N} \,=\,\Xmin{\eta}{\ss i 1 \hfill \ldots
\hfill N \\ \ss j k_1\hfill \ldots \hfill k_N}\ .$$ Eq.\eqref{detv}
follows by \eqref{vXeta}. \hfill$\square$

\section{Determinantal relations and combinatorial
$\theta$ identities}\label{secthetas}

In this section, after recalling basic definitions about theta
functions and fixing the necessary notation, we will derive some
identities, such as the inverse of the Brill-Noether matrix
expressed in terms of theta functions, Proposition \ref{inverr},
and, applying the results of the previous sections, will express the
relations among elements of $H^0(K_C^2)$
or, equivalently, the generators of the ideal of the canonical curve, in
terms of combinatorial theta identities. We will also introduce a
new basis for the holomorphic sections of $K_C$ expressed in terms
of Szeg\"o kernels. Such a basis satisfies properties that will lead
to an expression of the determinant of the Brill-Noether matrix in
terms of theta functions without using the problematic $\sigma$ section by Klein and Fay. Such a basis is also
the key point leading to Theorem \ref{szeghimatr}, whose proof, together with Theorems \ref{main},  \ref{ththetarel} and Corollary
\ref{threla}, is given in this section.
In particular, the dependence of the coefficients of the
quadrics and, possibly, cubics, that generate the ideal of a
canonical curve, on the points entering in Petri's construction, is
shown explicitly in terms of Riemann theta functions.

\vskip 6pt

Let $\Hh_g:=\{Z\in M_g(\CC)\mid \tp Z=Z,\im Z>0\}$, be the Siegel
upper half-space, i.e. the space of $g\times g$ complex symmetric
matrices with positive definite imaginary part, and define the usual
action of the symplectic group $\Gamma_g:={\rm Sp}(2g,\ZZ)$ on
$\Hh_g$ by
$Z\mapsto (AZ+B)(CZ+D)^{-1}$, $\Big(\begin{matrix}A & B\\ C &D\end{matrix}\Big)\in
\Gamma_g$.
Set ${\cal A}_Z:=\CC^g/L_Z$, $L_Z:={\ZZ}^g +Z{\ZZ}^g$, $Z\in \Hh_g$
and consider the theta function with characteristics
$$\theta \left[^a_b\right]\left(z,Z\right):= \sum_{k\in
{\ZZ}^g}e^{\pi i \tp{(k+a)}Z(k+a)+ 2\pi i \tp{(k+a)} (z+b)} \ ,
$$
where $z\in {\cal A}_Z$, $a,b\in{\RR}^g$. It has the
quasi-periodicity properties
$$\theta \left[^a_b\right]\left(z+n+Zm,Z\right)=
e^{-\pi i \tp{m}Zm-2\pi i\tp{m}z+2\pi i(\tp{a}n-\tp{b}m)}\theta
\left[^a_b\right]\left(z,Z\right) \ ,
$$
$m,n\in\ZZ^g$.

Denote by $\Theta\subset {\cal A}_Z$ the divisor of
$\theta(z,Z):=\theta\left[^0_0\right](z,Z)$ and by
$\Theta_{s}\subset \Theta$ its singular locus, where both $\theta$ and its gradient
vanish.

\subsection{Riemann theta functions and prime form}\label{thetaprime}

A different choice of the symplectic basis of $H_1(C,\ZZ)$
corresponds to a $\Gamma_g$ transformation
$$\left(\begin{matrix}\alpha\cr
\beta\end{matrix}\right)\mapsto \left(\begin{matrix}\tilde\alpha\cr
\tilde\beta\end{matrix}\right)=\left(\begin{matrix}D & C \cr B & A\end{matrix}\right)
\left(\begin{matrix}\alpha\cr \beta\end{matrix}\right)\ ,\qquad\qquad \left(\begin{matrix}A
& B \cr C & D\end{matrix}\right)\in \Gamma_g \ ,
$$
\begin{equation}\label{modull}\tau \mapsto
\tau'=(A\tau+B)(C\tau+D)^{-1}\ .\end{equation}
If $\delta',\delta''\in\{0,1/2\}^g$, then $\theta
\left[\delta\right]\left(z,\tau\right):=\theta
\left[^{\delta'}_{\delta''}\right]\left(z,\tau\right)$ has definite
parity in $z$ $\theta \left[\delta\right]
\left(-z,\tau\right)=e(\delta) \theta \left[\delta\right]
\left(z,\tau\right)$, where $e(\delta):=e^{4\pi i\!\tp{\delta'}
\delta''}$. There are $2^{2g}$ different characteristics of definite
parity, $2^{g-1}(2^g+1)$
even and $2^{g-1}(2^g-1)$ odd. By Abel Theorem each one of such characteristics determines
the divisor class of a spin bundle $L_\delta\simeq K^{1\over2}_C$,
so that we may call them spin structures. The dimension $h^0(L_\delta)$ of the space of global sections of $L_\delta$ has the same parity as the corresponding theta characteristic and equals the degree of the zero of $\theta \left[\delta\right]
\left(z,\tau\right)$ at $z=0$. For a generic curve, $h^0(L_\delta)$ is equal to zero or one; when this is not the case the spin structure (and the corresponding theta characteristic) is called singular. The curves admitting singular spin structures form a sublocus of codimension one in the moduli space.

Denote by $C_n=\Sym^n(C)$, $n\in\ZZ$, the space of effective divisors of degree $n$ on $C$. Let $J_n(C)$
be the principal homogeneous space of linear equivalence classes of divisors of degree $n$ on $C$ and
$J(C):={\CC}^g/L_\tau$, $L_\tau:={\ZZ}^g +\tau {\ZZ}^g$ the Jacobian of $C$.
Choose an arbitrary point $p_0\in C$ and denote by $A:C\to J(C)$,
$A(p):=(A_1(p),\ldots,A_g(p))$, $A_i(p):=\int_{p_0}^p\omega_i$,
$p\in C$, the Abel-Jacobi map.
$J(C)$ is identified with $J_0(C)$: each point of $J_0(C)$ can be expressed as $D_2-D_1$ with $D_1$ and $D_2$ effective divisors of the same degree,
this corresponds to $A(D_1-D_2)\in J(C)$, where $A(\sum_i n_i p_i):=\sum_in_iA(p_i)$, $p_i\in C$, $n_i\in\ZZ$.
Note that all the maps $C^g\to C_g\to J(C)$ are surjective.
 Let $W_n$ be
the closure of the image $A(C_n)$ in  $J_0(C)$. Denote by
$W_n^r\subseteq W_n$, $r\ge 0$, the image of the divisors $\d\in
C_n$ with $h^0(\O(\d))-1\ge r$, and observe that $W_n^0= W_n$.
Let $A\subset
C^g:=C\times\ldots\times C$ be the locus of $g$-tuples of points
$(p_1,\ldots,p_g)$ such that $\det\phi^1_i(p_j)=0$. Note that $A$
is independent of the choice of the basis of $H^0(K_C)$; in fact,
$A$ is the union of the hyperplanes $p_i=p_j$, $i\neq j\in I_g$,
and of the set of points such that $p_1+\ldots+p_g$ is a special
divisor, that is such that $\sum_{i=1}^gA(p_i)\in W_g^1$.
Consider the vector of Riemann constants
$$K^p_i:={1\over 2}+{1\over
2}\tau_{ii}-\sum_{j\neq
i}^g\oint_{\alpha_j}\omega_j\int_{p}^x\omega_i \ ,$$
$i\in I_g$, for any
$p\in C$. For any $p\in C$ define the formal sum $$\Delta:=(g-1)p-K^p=
(-{1\over 2}-{1\over
2}\tau_{ii}+\sum_{j\neq
i}^g\oint_{\alpha_j}\omega_j\int^x_{\,\cdot}\omega_i)_{i=1,\ldots,g} \ ,$$
so that, for any divisor $\xi$ of degree $g-1$ in $C$, $\xi-\Delta$ is the point in $\CC^g$
 given by $\int_{(g-1)p}^\xi \omega+K^P$. Under the projection $\CC^g\to J_0(C)$, $\Delta$ becomes a distinguished point in
 $J_{g-1}(C)$ depending only on the homological class of the marking. Furthemore,  $2\Delta =
K_C$. See \cite{jfayy}
for further details.
We will consider $\theta(\d+e):=\theta
\left[^0_0\right](A(\d)+e,\tau)$, for all $e\in J_0(C)$, evaluated
at some $0$-degree divisor $\d$ of $C$. We will also use the
notation
$$\deltadiv(\d):=\theta(A(\d-n\Delta))\ , $$
$$\theta_{n\Delta,i}(\d):=\partial_{z_i}\theta(z)_{|z=\int_{n(g-1)p}^\d\omega+nK^p} \ , $$
for each divisor $\d$ of degree $n(g-1)$, $n\in\ZZ$.
According to the Riemann Vanishing Theorem, for any $p\in C$ and
$e\in J_0(C)$

\smallskip
\begin{itemize}
\item[\it i.] if $\theta(e)\ne 0$, then the divisor $\d$ of $\theta(x-p-e)$
in $C$ is effective of degree $g$, with index of specialty $i(\d):=h^0(K_C\otimes \O(-\d))=0$
and $e=A(\d-p)-\Delta$;

\smallskip

\item[\it ii.] if $\theta(e)=0$, then for some $\zeta\in C_{g-1}$,
$e=A(\zeta)-\Delta$.

\end{itemize}

\vskip 6pt

\noindent Let $\nu$ be a non-singular odd characteristic. Since $\nu$ corresponds to an effective square root of the canonical divisor, it follows that
$$h_\nu(p)^2:=\sum_{1}^g\omega_i(p)
\partial_{z_i}\theta\left[\nu\right](z)_{|_{z=0}} \ ,$$ is the square of a holomorphic section $h_\nu(p)$ of the spin bundle $L_\nu$.
Therefore, the divisor of $h^2_\nu(p)$ is composed by $g-1$
double zeros  for all $p\in C$.
The prime form (see, for example, \cite{jfayy}, section II) $$
E(z,w):={\theta\left[\nu\right](w-z,\tau)\over
h_{\nu}(z)h_{\nu}(w)}\ ,
$$ is a holomorphic section of a line bundle
on $C\times C$, satisfying the multi-valuedness properties $$
E(z+\tp\alpha n+\tp\beta m,w)=\chi e^{-\pi i \tp m \tau m -2\pi i
\tp m A(z-w)}E(z,w) \ ,
$$
where $\chi:=e^{2\pi i(\tp \nu' n- \tp \nu'' m)}\in\{-1,+1\}$,
$m,n\in\ZZ^g$ and such that $E(z,w)=-E(w,z)$.
In particular, it only vanishes if $z=w$, and if $t$ is a local coordinate at $z\in
C$ such that $h_\nu=dt$, then
$$
E(z,w)={t(w)-t(z)\over\sqrt{dt(w)}\sqrt{dt(z)}}(1+{\cal
O}((t(w)-t(z))^2)) \ .
$$

Let $\delta:=\left[^{\delta'}_{\delta''}\right]\in \ZZ_2^{2g}$ be an
even theta characteristic on $J_0(C)$ for some curve $C$ of genus
$g$ and denote by $L_\delta$ the corresponding spin bundle. The Szeg\"o kernel $S_\delta$ is the unique
section of $L_\delta\boxtimes L_\delta$ on $C\times C$, with a
simple pole with residue $1$ at the diagonal and holomorphic
otherwise (see \cite{jfayy}, pg 26). It can be expressed in terms of the prime form as
\begin{equation}\label{szegodef}S_\delta(p,q)={\theta[\delta](p-q)\over\theta[\delta](0)E(p,q)}\ . \end{equation}
The Szeg\"o kernel satisfies the  Fay
identity for more points \cite{jfayy}, also known as Gunning secant formula \cite{gunningtwo} (see also \cite{BL} and references therein)
\begin{equation}\label{Fay}{\theta[\delta](\sum_{1}^m(x_i-y_i))\prod_{i<j}E(x_i,x_j)E(y_i,y_j)\over
\theta[\delta](0)\prod_{i,j}E(x_i,y_j)}=(-1)^{{m(m-1)\over
2}}\det\nolimits_{ij}{S_\delta(x_i,y_j)}\ ,\end{equation} for any non singular
even spin structure $\delta$ and arbitrary points
$x_1,\ldots,x_m$, $y_1,\ldots,y_m\in C$. The case $n=2$ corresponds to the Fay's trisecant identity
which admits other basic generalizations introduced by Gunning in \cite{gunningtwo}.

Denote by $\omega_{a-b}(x)$ $a,b,x\in C$, the meromorphic
$1$-differential of the second kind with simple poles on $a$ and $b$
with residue $+1$ and $-1$, respectively, holomorphic on
$C\setminus\{a,b\}$. Choose a base point $p_0\in C$ and a set of
generators $\alpha_1,\ldots,\alpha_g,\beta_1,\ldots,\beta_g$ of the
first homotopy group $\pi(C,p_0)$, satisfying the unique relation
$\prod_{i=1}^g\alpha_i\beta_i\alpha_i^{-1}\beta_i^{-1}=1$, and
consider the associated dissection $\tilde C$ of $C$ along the
representatives of such generators. Then, the differential
$\omega_{a-b}$, for all $a,b$ in the interior of $\tilde C$, is
completely determined by imposing the conditions
\begin{equation}\label{ominta}\oint_{\alpha_i}\omega_{a-b}=0\ , \end{equation} $i\in I_g$.
Furthermore,  the following properties
\begin{equation}\label{omintb}\oint_{\beta_i}\omega_{a-b}=\int_a^b\omega_i\ ,\end{equation} hold
for all $i\in I_g$, where the integrations are on $\tilde C$ \cite{jfayy}.

We will also consider the multi-valued $g/2$-differential
$\sigma(z)$ on $C$ with empty divisor \cite{FayMAM}, satisfying the property
$$
\sigma(z+\tp\alpha n+\tp\beta m)=\chi^{-g}e^{\pi i(g-1)\tp m \tau
m+2\pi i \tp m K^z}\sigma(z)\ .
$$
Such conditions fix $\sigma(z)$ only up to a factor independent of
$z$; the precise definition, to which we will refer, can be given,
following \cite{FayMAM},
on the universal covering of $C$  (see also \cite{jfayy}). Furthermore \cite{FayMAM}
$$
\sigma(z,w):={\sigma(z)\over\sigma(w)}={\deltadiv(\sum_{1}^gx_i-z)\over
\deltadiv(\sum_{1}^gx_i-w)}\prod_{i=1}^g{E(x_i,w)\over E(x_i,z)} \ ,$$
for all $z,w,x_1,\ldots,x_g\in C$, which follows by observing that the RHS multiplied by $\sigma(w)/\sigma(z)$ is a meromorphic function with empty divisor.

\subsection{Determinants and theta functions on algebraic curves}

For all $y,p_1,\ldots,p_g\in C$ set  \cite{jfayy} \begin{equation}\label{essesimm}S(p_1,\ldots,p_g):={\deltadiv(\sum_1^gp_i-y)\over
\sigma(y)\prod_{1}^gE(y,p_i)}\ . \end{equation}

\noindent The following properties of $S$ are consequence of standard arguments based, essentially, on the Riemann Vanishing Theorem.

\begin{lem}\label{lemessesimm}For all $p_1,\ldots,p_g\in C$, $S(p_1,\ldots,p_g)$ is independent of $y$.
For each fixed $\d\in C_{g-1}$, consider the map $\pi_\d:C\to C_g$,
$\pi_\d(p):=p+\d$. The pull-back $\pi_\d^*S$ vanishes identically if
and only if $\d$ is a special divisor; if $\d$ is not special, then
$\pi_\d^*S$ is the unique (up to normalization) holomorphic
$1/2$-differential, with a
fixed spin structure,
such that $[(\pi_\d^*S)+\d]$ is the canonical
divisor class.\end{lem}

\vskip 6pt

\begin{remark}\label{questor}
The fact that $S(p_1,\ldots,p_g)=0$ if and only if $A(p_1+\ldots+p_g)\in W^1_g$,
implies that $S$ can be defined unambiguously as a function on $W_g$. More precisely, although $S$ apparently depends on the precise point in $C_g$,
it just vanishes on the locus where the Jacobian map fails to be one to one.
\end{remark}

\vskip 6pt

\noindent
Checking the divisors and noticing that the LHS are ratios of the same sections, one may easily see that the following proposition holds.

\vskip 6pt

\begin{prop}\label{thdettheta}For each $n\in\ZZ_{>0}$, let $\phi^n:=\{\phi_i^n\}_{i\in
I_{N_n}}$ be an arbitrary basis of $H^0(K_C^n)$ and $p_1,\ldots,p_{N_n}$ arbitrary points in $C$. For $n\geq2$, there are constants
$\kenne[\phi^n]$, depending only on the marking of $C$ and on
$\{\phi_i^n\}_{i\in I_{N_n}}$ but independent of $p_1,\ldots,p_{N_n}$, such that
  \begin{equation}\label{dettheta}\kuno[\phi^1]=
{\det \phi_i^1(p_j)\over
S(p_1,\ldots,p_g)\prod_1^g\sigma(p_i)
\prod_{i<j}^gE(p_i,p_j) }\ ,\end{equation} and
\begin{equation}\label{detthetaii}\kenne[\phi^n]={\det \phi_i^{n}(p_j)\over
\theta_\Delta\bigl(\sum_{1}^{N_n}
p_i\bigr)\prod_{1}^{N_n}\sigma(p_i)^{2n-1}\prod_{i<j}^{N_n}
E(p_i,p_j)}\ .\end{equation} \end{prop}

\vskip 6pt

\begin{remark}\label{usefulassaie} In the case $\phi_i=\omega_i$, $i\in I_g$, Eq.\eqref{dettheta} corresponds to Corollary 2.17 of \cite{jfayy} (see also
Corollary 1.4 in \cite{FayMAM}). Note that Proposition \ref{thdettheta} also
implies that $S(p_1,\ldots,p_g)$ does not depend on $y$. By taking the
limit $y\to z:=p_i$ in $S(p_1,\ldots,p_g)$, it follows by \eqref{essesimm}
and \eqref{dettheta}, that \begin{equation}\label{ildet}\det \eta(z,p_1,\ldots,\check
p_i,\ldots,p_g) = \kuno[\eta]\sum_{l=1}^g\theta_{\Delta,l}(\a_i)
\omega_l(z)\prod_{{j,k\neq i \atop j<k}}E(p_j,p_k)\prod_{j\neq
i}\sigma(p_j)\ ,\end{equation} for all
$p_1,\ldots,p_{i-1},z,p_{i+1},\ldots,p_g\in C$, where
$\theta_{\Delta,i}(e):=\partial_{z_i}\theta_\Delta(z)_{|z=e}$,
$e\in\CC^g$ and $\a_i:=\sum_{j\neq i}p_j$, $i\in I_g$. Furthermore,
$$\theta(K^w+w-z)={\sigma(z)E(z,w)^g\over\kappa[\omega]\sigma(w)^g}W[\omega](w)\ ,
$$ for all $w,z\in C$, with $W[\omega](z)$ the Wronskian of
$\{\omega_i\}_{i\in I_g}$ at $z$. Note that by \eqref{ildet} the condition
$\det \eta_i(p_j)\neq 0$ implies
$$\sum_j\theta_{\Delta,j}(\a_i)\omega_j(p_i)\neq 0\ ,
$$ for
all $i\in I_g$. On the other hand, the LHS is in
$H^0(K_C)$, so that it may vanish either at the $2g-2$ points, or
identically if $A(\a_i)-\Delta\in\Theta_{s}$, $i\in I_g$.\end{remark}

\noindent Remarkably, the inverse of the Brill-Noether matrix
$[\omega]_{ij}:=\omega_i(p_j)$ admits an expression in terms of
theta functions. A key feature of this result, compared with analogous expressions provided, for example, by Fay \cite{jfayy}, is that such relationship does not involve the definition of trascendental functions such as the Klein-Fay section $\sigma$, nor the definition of moduli-dependent overall factors analogous to $\kuno[\phi^i]$ in Eq.\eqref{dettheta}. Thanks to the definition of the basis $\{\sigma_i\}_{i\in I_g}$, the derivation is surprisingly simple. This result
will be used to prove Theorem \ref{thleX}.

\vskip 6pt

\begin{prop}\label{inverr}Let $p_1,\ldots,p_g\in C$ be points such that
$\det \eta_i(p_j)\neq 0$. We have
\begin{equation}\label{leomega}[\omega]^{-1}_{ij}=\oint_{\alpha_j}\sigma_i={\theta_{\Delta,j}
\left(\a_i\right)\over\sum_k\theta_{\Delta,k}\left(\a_i\right)\omega_k(p_i)}\
,\end{equation} $i,j\in I_g$, so that
\begin{equation}\label{siggmas}\sigma_i(z)
=\sigma(z,p_i){\deltadiv(\a+z-y-p_i)\over\deltadiv(\a-y)E(z,p_i)}{E(y,p_i)\over
E(y,z)} \prod_{1}^gE(z,p_j)\prod_{j\neq i}{1\over E(p_i,p_j)} \ ,\end{equation}
and \begin{equation}\label{kkkk}\kuno[\sigma]={\sigma(y)\prod_1^gE(y,p_i)\over
\deltadiv(\a-y)\prod_{i<j}^gE(p_i,p_j)\prod_1^g\sigma(p_k)} \ ,\end{equation} for
all $z,y,x_i,y_i\in C$, $i\in I_g$,
 where $\a:=\sum_1^gp_i$, $\a_i:=\a-p_i$, $i\in I_g$.\end{prop}

\vskip 6pt

\noindent {\it Proof.} The matrix
$[\phi^n]_{ij}:=\phi^n_i(p_j)$ is non-singular, so that
$\gamma^n_i=\sum_{j=1}^{N_n}[\phi^n]^{-1}_{ij}\phi^n_j$
and by \eqref{newbasis}
$\sigma_i=\sum_j[\omega]^{-1}_{ij}\omega_j$, and \eqref{leomega}
follows by \eqref{dettheta} and \eqref{ildet}. Eq.\eqref{siggmas}
follows by \eqref{dettheta} and
\eqref{kkkk} by $\det \sigma_i(p_j)=1$.
\hfill$\square$

\vskip 6pt

Given a set of points $p_1,\ldots,p_g$, set $\a_i:=\sum_1^gp_k-p_i$ and $\b:=\sum_3^gp_i$ and define
\begin{equation}\label{leX}X^\omega_{ij}:={\theta_{\Delta,\one_j}(\a_{\one_i})
\theta_{\Delta,\two_j}(\a_{\two_i})+\theta_{\Delta,\one_j}(\a_{\two_i})
\theta_{\Delta,\two_j}(\a_{\one_i})\over
(1+\delta_{\one_j\two_j})\sum_{l,m}\theta_{\Delta,l}(\a_{\one_i})
\theta_{\Delta,m}(\a_{\two_i})\omega_l(p_{\one_i})\omega_m(p_{\two_i})}\
,\end{equation} $i,j\in I_M$. (Recall that $(\one_i,\two_i)$ denotes the $i$-th element in the ordered set $\{(1,1),\ldots,(g,g),(1,2),\ldots,(1,g),\ldots,(g-1,g)\}$, see Section \ref{notablecombina}). Let us define the holomorphic quadratic differentials
\begin{equation}\label{vXww}v_i:=\sum_{j=1}^MX^\omega_{ji}\,\omega\omega_j\ ,\end{equation} $i\in
I_N$, with $X^\omega_{ji}$ as in Eq.\eqref{leX}.

\begin{theorem}\label{thleX}Fix $g-2$ distinct points $p_3,\ldots,p_g\in
C$ such that \begin{equation}\label{ipotesi}\{A(p+\b)-\Delta|p\in C\}\cap
\Theta_{s}=\emptyset\ ,\end{equation} $\b:=\sum_3^gp_i$.
Then, for any choice of $p_2$ distinct from $p_3,\ldots,p_g$ and $p_1\in C\setminus S$, where $S$ is a finite set of points containing $p_2,\ldots,p_g$, the quadratic differentials $v_1,\ldots,v_{N}$ given in \eqref{vXww} form a
 basis for $H^0(K_C^2)$ independent of the choice of the marking. Furthermore, the $M-N$
independent linear relations
$$\sum_{j=1}^MC^\omega_{ij}\omega\omega_j=0\ ,
$$
$N+1\le
i\le M$, hold, with
\begin{equation}\label{leCi}C^\omega_{ij}:=\sum_{k_1,\ldots,k_N=1}^M \Xmin{\omega}{\ss
1 \hfill \ldots \hfill N i \\ \ss k_1\hfill \ldots \hfill k_N
j\hfill}\; \R_v[\omega\omega_{k_1},\ldots,\omega\omega_{k_N}]\ . \end{equation}
$X^\omega_{ij}$ and $C^\omega_{ij}$ are the coefficients defined in
\eqref{XXXX} and \eqref{leC}, respectively, for $\eta_i=\omega_i$, $i\in I_g$. \end{theorem}

\vskip 6pt

\noindent {\sl Proof.} After noting that $h^0(K_C\otimes\O(-\b))\ge
2$, fix a pair of linearly independent elements $\sigma_1,\sigma_2$
of $H^0(K_C\otimes\O(-\b))$. Eq.\eqref{ipotesi} implies that
$h^0(K_C\otimes\O(-\b-p))=1$, for all $p\in C$, so that
$\{\sigma_1,\sigma_2\}$ is a basis of $H^0(K_C\otimes\O(-\b))$ and
$\supp((\sigma_1)-\b)\cap \supp((\sigma_2)-\b)$ is empty. Fix
$p_2\in C\setminus\{p_3,\ldots,p_g\}$; without loss of generality,
we can assume $\sigma_1\in H^0(K_C\otimes\O(-\b-p_2))$. Define the
finite set $S$ as the support of $(\sigma_1)$ or, equivalently, as
the union of $\{p_2,\ldots,p_g\}$ and the set of zeros of
$S(x,p_2,\ldots,p_g)$. Then, for each $p_1\in C\setminus S$, by choosing
$\sigma_2\in H^0(K_C\otimes\O(-\b-p_1))$, we obtain that $\sigma_1$
and $\sigma_2$ are linearly independent. Then $p_1,\ldots,p_g$
satisfy the conditions {\it i}) and {\it ii}) of Proposition \ref{thlev},
and $\{v_i\}_{i\in I_N}$, as defined in \eqref{lev}, is a basis of
$H^0(K_C^2)$. Eq.\eqref{leomega} implies that Eq.\eqref{leX} is equivalent to
\eqref{XXXX}, and the theorem follows by Theorem \ref{thcorol}. \hfill$\square$

\subsection{Combinatorial theta identities on the canonical curve}

Applying Lemmas \ref{thcombi} and \ref{thcombvi} to determinants of
symmetric products of sections of $H^0(K_C)$ leads to
combinatorial relations that, by \eqref{dettheta} and \eqref{detthetaii}, yield
non trivial identities among theta functions.

\noindent
The relation $\det\omega\omega(x_1,\ldots,x_M)=0$, $g\geq 4$, Lemma \ref{thcombi} and
Eq.\eqref{dettheta} imply the following proposition.

\vskip 6pt

\begin{prop}\label{thdetvan}The following identity
$$\sum_{s\in\perm_M}\sgn(s) \prod_{i=1}^{g+1}\det
\omega(x_{d^i(s)})=0\ , \qquad\quad g\geq4\ ,
$$
where
$x_i$, $i\in I_M$, are arbitrary points of $C$, holds. Furthermore,
it is equivalent to
$$
\sum_{s\in\perm_M}\sgn(s)\prod_{k=1}^{g+1}
{\deltadiv\bigl(\sum_{i=1}^gx_{d^k_i(s)}-y_{k,s})\prod_{i<j}^gE(x_{d_i^k(s)},x_{d_j^k(s)})
\over \prod_{i=1}^gE(y_{k,s},x_{d^k_i(s)})\sigma(y_{k,s})}=0\ ,
$$
for $g\geq4$, where $y_{k,s}$, $k\in I_{g+1}$,
$s\in \perm_M$, are arbitrary points of $C$.\end{prop}

\vskip 6pt

\noindent The basis $\{\sigma_i\}_{i\in I_g}$ of $H^0(K_C)$ adapted to the points $p_1,\ldots,p_g$, satisfies the conditions Eq.\eqref{combi} of Lemma \ref{thcombi}, for any $n<g$. Thus, we can apply Lemmas \ref{thcombi} and \ref{thcombvi} to the determinants of quadratic differentials of the form $\sigma_i\sigma_j$. This leads to an expression in terms of determinants of holomorphic $1$-differentials only.

\vskip 6pt

\begin{theorem}\label{thdetvp}Fix the points $p_1,\ldots,p_g\in C$,
and $\hat\sigma_i\in H^0(K_C)$, $i\in I_g$, in such a way that
$\hat\sigma_i(p_j)=0$, for all $i\neq j\in I_g$. Define $\hat v_i\in
H^0(K_C^2)$, $i\in I_N$, by
$$\hat v_i:=\psi(\spbase{\hat\sigma}_i)=
\hat\sigma_{\one_i}\hat\sigma_{\two_i}\ ,$$ and let $\{\eta_i\}_{i\in
I_g}$ be an arbitrary basis of $H^0(K_C)$. Then, the following
identity \begin{align}\label{dettv}&\det \hat
v(p_3,\ldots,p_g,x_1,\ldots,x_{2g-1})\Bigl({\det
\eta_i(p_j)\over\hat\sigma_1(p_1)\hat\sigma_2(p_2)}\Bigr)^{g+1}
\prod_{i=3}^g\hat\sigma_i(p_i)^{-4}\notag \\ &={(-1)\over
c_{g,2}}^{g+1}\sum_{s\in\perm_{2g-1}}\sgn(s)\det
\eta(x_{d^1(s)})\det \eta(x_{d^2(s)})\notag \\ & \cdot \prod_{i=3}^{g+1}\det
\eta(x_{d^i_1(s)},x_{d^i_2(s)},p_3,\ldots,p_g)\ ,\end{align} holds for all
$x_1,\ldots,x_{2g-1}\in C$, where, according to \eqref{glic},
$c_{g,2}=g!(g-1)!(2g-1)$.\end{theorem}

\vskip 6pt

\noindent {\sl Proof.} Assume that $p_1,\ldots,p_g$ satisfy the
hypotheses of Proposition \ref{thnewbasis}, so that
$\{\hat\sigma_i\}_{i\in I_g}$ is a basis of $H^0(K_C)$ adapted to $p_1,\ldots,p_g$. In this case, \eqref{dettv} follows by Lemma \ref{thcombi} for $n=2$ and by the identity $\det \hat\sigma_i(z_j)={\det\eta_i(z_j)\over \det
\eta_i(p_j)}\prod_{i=1}^g\hat\sigma_i(p_i)$. Since \eqref{dettv} holds for a dense set of $g$-tuples $p_1,\ldots,p_g$ in
$C^g$, we conclude by continuity arguments.
\hfill$\square$

\vskip 6pt

\begin{remark}\label{thdetvxxdx} If $\det \eta_i(p_j)\neq 0$, then Theorem
\ref{thdetvp} holds for $\hat\sigma_i=\sigma_i$, so that
$\hat\sigma_i(p_i)=1$, $i\in I_g$, and $\hat v_i= v_i$, $i\in I_N$.\end{remark}

\vskip 6pt

\begin{cor}\label{thdetvx}Let $\b:=\sum_{3}^gp_i$ be a fixed divisor of $C$ and define $\hat v_i$, $i\in I_N$, as in Theorem
\ref{thdetvp}. Then for all $x_1,\ldots,x_N\in C$
\begin{align}\label{detvi}&\det \hat v(x_1,\ldots,x_N)=F{\deltadiv\bigl(\sum_{1}^Nx_i\bigr)\prod_{i=2g}^N(\sigma(x_i)^3
\prod_{j=1}^{i-1}E(x_j,x_i))\over
\deltadiv\bigl(\sum_{1}^{2g-1}x_i+\b\bigr)\prod_{i=3}^{g}
\prod_{j=1}^{2g-1}E(p_i,x_j)}\prod_{i=1}^{2g-1}\sigma(x_i)^2\notag\\
\cdot&\sum_{s\in\perm_{2g-1}}\sgn(s)S(x_{s_1},\ldots,x_{s_g})
S(x_{s_g},\ldots,x_{s_{2g-1}})
\prod_{^{i,j=1}_{i<j}}^gE(x_{s_i},x_{s_j})\prod_{^{i,j=g}_{i<j}}^{2g-1}
E(x_{s_i},x_{s_j}) \notag\\\cdot
&\prod_{k=1}^{g-1}\Bigl(S(x_{s_k},x_{s_{k+g}},p_3,\ldots,p_g)
E(x_{s_k},x_{s_{k+g}})\prod_{i=3}^gE(x_{s_k},p_i)E(x_{s_{k+g}},p_i)\Bigr)\
,\end{align} where
\begin{align}F &:= -{1\over
c_{g,2}}\Bigl(
{\hat\sigma_1(p_1)\hat\sigma_2(p_2)\over
S(p_1,\ldots,p_g)\sigma(p_1)\sigma(p_2)E(p_1,p_2)}\Bigr)^{g+1}\notag \\
&\cdot\prod_{i=3}^g{\hat\sigma_i(p_i)^4\over \sigma(p_i)^5
(E(p_1,p_i)E(p_2,p_i))^{g+1}\prod_{j>i}^gE(p_i,p_j)^3} \ . \notag \end{align}\end{cor}

\vskip 6pt

\noindent {\sl Proof.} Apply Eq.\eqref{dettv} to
$$\det \hat
v(x_1,\ldots,x_N)={\det \rho(x_1,\ldots,x_N)\det \hat
v(p_3,\ldots,p_g,x_1,\ldots,x_{2g-1})\over \det
\rho(p_3,\ldots,p_g,x_1,\ldots,x_{2g-1})}\ ,$$ with $\{\rho_i\}_{i\in
I_N}$ an arbitrary basis of $H^0(K_C^2)$. Eq.\eqref{detvi} then follows by
\eqref{dettheta} and \eqref{detthetaii}. \hfill$\square$

\begin{theorem}\label{thlaH}Fix $p_1,\ldots,p_g\in C$.
\begin{align}\label{ilrappo} & H(p_1,\ldots,p_g):= \notag \\
& {S(p_1,\ldots,p_g)^{5g-7}E(p_1,p_2)^{g+1}\over
\deltadiv\bigl(\b+\sum_{1}^{2g-1}x_i\bigr)\prod_{i=1}^{2g-1}
\sigma(x_i)}
\prod_{i=3}^g{E(p_1,p_i)^4E(p_2,p_i)^4\prod_{j>i}^gE(p_i,p_j)^5\over
\sigma(p_i)}\notag \\ &\cdot \sum_{s\in\perm_{2g-1}}
{S(x_{s_1},\ldots,x_{s_g})
S(x_{s_g},\ldots,x_{s_{2g-1}})\over
\prod_{i=3}^gE(x_{s_g},p_i)}
\prod_{i=1}^{g-1}{S(x_{s_i},x_{s_{i+g}},p_3,\ldots,p_g)\over \prod_{_{j\neq
i}^{j=1}}^{\ss g-1}E(x_{s_i}, x_{s_{j+g}})} \ ,\end{align} is independent of
the points $x_1,\ldots,x_{2g-1}\in C$. Furthermore, the set
$\{v_i\}_{i\in I_N}$, defined as in {\rm \eqref{lev}}, is a basis of
$H^0(K_C^2)$ if and only if $H(p_1,\ldots,p_g)\neq 0$.\end{theorem}

\vskip 6pt

\noindent {\sl Proof.} Consider
$$\hat\sigma_i(z):=A_i^{-1}\sigma(z)S(p_1,\ldots,\check p_i,\ldots,p_g,z)
\prod_{^{j=1}_{j\neq
i}}^gE(z,p_j)=A_i^{-1}\sum_{j=1}^g\theta_{\Delta,j}(\a_i)\omega_j(z)\
,$$ $i\in I_g$, $\a_i:=\a-p_i$, with $A_1,\ldots,A_g$ non-vanishing
constants. Such elements of $H^0(K_C)$ are independent of $y$ and if
the points $p_1,\ldots,p_g$ satisfy the hypotheses of Proposition
\ref{thnewbasis}, then, recalling that
$\sigma_i=\sum_j[\omega]_{ij}^{-1}\omega_j$, with
$\{\sigma_i\}_{i\in I_g}$ defined in \eqref{newbasis}, it follows by
\eqref{leomega} that the basis $\sigma_1,\ldots,\sigma_g$ of $H^0(K_C)$
coincides, up to a non-singular diagonal transformation, with
$\hat\sigma_1,\ldots,\hat\sigma_g$. Let $\{\rho_i\}_{i\in I_N}$ be
an arbitrary basis of $H^0(K_C^2)$. By \eqref{detthetaii} the following
identity
\begin{align}&\det \rho(p_3,\ldots,p_g,x_1,\ldots,x_{2g-1}) = \notag \\ &
\kdue[\rho] \sgn(s)\deltadiv\Bigl(\smsum_{\ss 1}^{\ss
2g-1}x_i+\b\Bigr)\prod_{_{i<j}^{i,j=1}}^{2g-1}E(x_{s_i},x_{s_j})
 \prod_{i=1}^{2g-1}\sigma(x_i)^3\prod_{i=3}^g\sigma(p_i)^3 \notag \\ &
\cdot\prod_{_{i<j}^{i,j=3}}^gE(p_i,p_j)\prod_{i=3}^g\prod_{j=1}^{2g-1}E(p_i,x_j)
\ ,\notag \end{align}
holds for all $s\in\perm_{2g-1}$. Together with Eq.\eqref{detvi} and the
above expression for $\hat\sigma_i$, it implies that
\begin{equation}\label{ratioo}H(p_1,\ldots,p_g)=\kappa[\rho]c_{g,2}
(A_1A_2)^{g+1}\prod_{i=3}^gA_i^4 {\det \hat
v(p_3,\ldots,p_g,x_1,\ldots,x_{2g-1})\over \det
\rho(p_3,\ldots,p_g,x_1,\ldots,x_{2g-1})}\ .\end{equation} Hence,
$H(p_1,\ldots,p_g)$ is independent of $x_1,\ldots,x_{2g-1}$, and
$H(p_1,\ldots,p_g)\neq0$ if and only if $\{\hat v_i\}_{i\in I_N}$ is
a basis of $H^0(K_C^2)$. On the other hand, the vector $(\hat
v_1,\ldots,\hat v_N)$ corresponds, up to a non-singular diagonal
transformation, to $(v_1,\ldots,v_N)$, with $v_i$, $i\in I_N$,
defined in \eqref{lev}. \hfill$\square$

\vskip 6pt

\begin{remark}\label{gfgfg}By \eqref{ratioo}
$$\kappa[\hat v]={H(p_1,\ldots,p_g)\over c_{g,2}
(A_1A_2)^{g+1}\prod_{i=3}^gA_i^4}\ .$$ Furthermore, if
$(p_1,\ldots,p_g)\notin A$, then one can choose
$$A_i=\sigma(p_i)S(p_1,\ldots,p_g)\prod_{^{j=1}_{j\neq
i}}^gE(p_i,p_j)=\sum_{j=1}^g\theta_{\Delta,j}(\a_{i})\omega_j(p_i) \ ,$$
to obtain $\hat\sigma_i=\sigma_i$, $i\in I_g$, and
\begin{align}\label{kkl}\kappa[v] & ={H(p_1,\ldots,p_g)\over
c_{g,2}S(p_1,\ldots,p_g)^{6g-6}\prod_{i=1}^2\bigl(\sigma(p_i)\prod_{^{j=1}_{j\neq
i}}^gE(p_i,p_j)\bigr)^{g+1}} \notag \\
&\cdot {1\over \prod_{i=3}^g\bigl(\sigma(p_i)\prod_{^{j=1}_{j\neq
i}}^gE(p_i,p_j)\bigr)^4}\ .\end{align}\end{remark}

\vskip 6pt

\begin{cor}\label{corollo}\begin{align}\label{iltrap}K(p_3,\ldots,p_g)&:={1\over
\deltadiv\bigl(\b+\sum_{1}^{2g-1}x_i\bigr)\prod_{1}^{2g-1}\sigma(x_i)\prod_{i=3}^g\sigma(p_i)
}\notag \\ &\cdot \sum_{s\in\perm_{2g-1}}
{S(x_{s_1},\ldots,x_{s_g})
S(x_{s_g},\ldots,x_{s_{2g-1}}) \over
\prod_{i=3}^gE(x_{s_g},p_i)} \notag \\
&\cdot \prod_{i=1}^{g-1}{
S(x_{s_i},x_{s_{i+g}},p_3,\ldots,p_g) \over \prod_{_{j\neq i}^{j=1}}^{\ss
g-1}E(x_{s_i}, x_{s_{j+g}})} \ , \end{align} is independent of
$x_1,\ldots,x_{2g-1}\in C$. Furthermore, for any $p_1,\ldots,p_g\in
C$ such that $\det \eta_i(p_j)\neq 0$, the set $\{v_i\}_{i\in I_N}$,
defined in {\rm \eqref{lev}}, is a basis of $H^0(K_C^2)$ if and only if
$K(p_1,\ldots,p_g)\neq 0$.\end{cor}

\vskip 6pt

\noindent {\sl Proof.} The ratio \begin{equation}\label{baaa}{H\over
K}=S(p_1,\ldots,p_g)^{5g-7}E(p_1,p_2)^{g+1}\prod_{i=3}^g(E(p_1,p_i)E(p_2,p_i))^4
\prod_{i<j=3}^gE(p_i,p_j)^5\ ,\end{equation} is independent of
$x_1,\ldots,x_{2g-1}$, so that the first statement follows by
Theorem \ref{thlaH} or, equivalently, noticing that by
Eqs.\eqref{kkkk}, \eqref{ilrappo}, \eqref{kkl}, \eqref{baaa}
$$K(p_3,\ldots,p_g)=(-1)^{g+1}c_{g,2}{\kdue[v]\over\kuno[\sigma]^{g+1}}\prod_{_{i<j}^{i,j=3}}^gE(p_i,p_j)^{2-g}
\prod_{i=3}^g\sigma(p_i)^{3-g}\ . $$ By \eqref{dettheta}, \eqref{essesimm} and \eqref{baaa}
the condition $\det \eta_i(p_j)\neq0$ implies $H/K\neq 0$. In this
case $K\neq 0$ if and only if $H\neq 0$, and the corollary follows
by Theorem \ref{thlaH}.\hfill$\square$

\subsection{Determinant of the Brill-Noether matrix in terms of theta
functions}

\begin{theorem}\label{ththetarel}Fix $g-2$ distinct points $p_3,\ldots,p_g\in C$ in such
a way that $K(p_3,\ldots,p_g)\neq 0$. The following $(g-2)(g-3)/2$
independent relations
\begin{align}\label{thetarel} V_{i_1i_2}(p_3,\ldots,p_g,&x_1,\ldots,x_{2g}):=\notag\\
\sum_{s\in\perm_{2g}}\sgn(s) \Biggl\{&\prod_{k=1}^2{S(\hat
x_{k},\hat x_{{g+k}},\hat x_{{2g}},p_3,\ldots,\check p_{i_k},\ldots,p_g)E(\hat x_{k},\hat
x_{{2g}}) E(\hat x_{{k+g}},\hat x_{{2g}})\over E(\hat
x_{k},p_{i_k})E(\hat x_{{k+g}},p_{i_k})E(\hat x_{{2g}},p_{i_k})}\notag\\
&\cdot\prod_{k=1}^{g-1}\bigl(E(\hat x_{k},\hat x_{{k+g}})
\prod_{j=3}^gE(\hat x_{k},p_j) E(\hat x_{{k+g}},p_j)\bigr) \notag \\
&\cdot S(\hat x_{1},\ldots,\hat x_g)
\prod_{_{k<j}^{k,j=1}}^gE(\hat x_{k},\hat x_{j})
S(\hat
x_{g},\ldots,\hat x_{2g-1})\prod_{_{k<j}^{k,j=g}}^{2g-1}E(\hat x_{k},\hat x_{j}) \notag\\
&\cdot\prod_{k=3}^{g-1}S(\hat x_{k},\hat
x_{{k+g}},p_3,\ldots,p_g)\prod_{j=3}^gE(\hat x_{{2g}},p_j)^2\Biggr\}=0 \ ,
\end{align}$3\le i_1<i_2\le g$, where $\hat x_i:=x_{s(i)}$, $i\in I_{2g}$,
hold for all $x_i\in C$, $i\in I_{2g}$.\end{theorem}

\vskip 6pt

\noindent {\sl Proof.} $V_{ij}(p_3,\ldots,p_g,x_1,\ldots,x_{2g})$
is obtained by expressing the determinants of the Brill-Noether matrices in the identities of Theorem \ref{main} in terms of theta functions
by means of \eqref{dettheta}. \hfill$\square$

\vskip 6pt

Note that $V_{ii}\neq 0$ for $i=3,\ldots,g$, since the expression on
the LHS is proportional to a determinant of $2g$ linearly
independent holomorphic sections of $H^0(K_C^2)$, evaluated at
general points $x_k\in C$, $k\in I_{2g}$.

By a limiting procedure we are able to derive the original Petri's relations, but
now written in terms of the basis $\{\omega_i\}_{i\in I_g}$ and with the coefficients explicitly expressed in terms of theta functions.
In particular, by \eqref{thetarel} one can easily write each term in the sum over permutations as a linear
combination of products $\omega_i(x_{2g})\omega_j(x_{2g})$, $i,j=1,\ldots,g$, whose coefficients are independent of $x_{2g}$.
Furthermore, remarkably, the dependence on $p_1,p_2$ and $x_1,\ldots,x_{2g-1}$ disappears by simply multiplying $V_{i_1i_2}$ by a suitable factor. This is the main idea behind the next corollary.
Before formulating it we need to recall few facts about the Kodaira-Spencer map.

Let $\T_g$ denote the Torelli space of smooth curves of
genus $g$ with a marking
and consider the embedding
$\pi:\T_g\rightarrow\Hh_g$ which associates to each marked $C$ its
Riemann period matrix $\tau$. Let $\hat\T_g\subset\T_g$ be the space
of marked canonical curves. For each $C\in\hat\T_g$ the fiber
$T^*_C\hat\T_g$ of the cotangent to $\T_g$ at the point
corresponding to $C$ is spanned by the pull-backs
$d\tau_{ij}:=\pi^*dZ_{ij}$, $i,j\in I_g$, of the elements
$dZ_{ij}\in T^*_{\pi(C)}\Hh_g$.

The two-fold products $\omega_i\omega_j$ are directly
related, via the Kodaira-Spencer map, to the elements $d\tau_{ij}$
on $T^*_C\hat\T_g$. More precisely,
consider the Kodaira-Spencer map $k$ identifying
$H^0(K_C^2)$ with the fiber of the cotangent bundle of
moduli space $\M_g$ of the smooth Riemann surfaces, isomorphic to $\T_g/\Gamma_g$, at the
point representing $C$. Next, consider a Beltrami differential
$\mu\in\Gamma(\bar K_C\otimes K_C^{-1})$
and recall that it defines a tangent vector at $C$ of $\T_g$.
The derivative of the period map $\tau_{ij}:\T_g\to \CC$ at $C$ in the direction
of $\mu$ is given by Rauch's formula \cite{FayMAM}
$d_C\tau_{ij}(\mu)=\int_C\mu\omega_i\omega_j$, so that $\forall j,k\in I_g$,
\begin{equation}\label{kodairaspecncer} k(\omega_j\omega_k)={1\over2\pi i}d\tau_{jk} \ .\end{equation}

\noindent Recall that the elements in the set
$\{\omega\omega_i\}_{i\in I_M}$ are $\omega_j\omega_k$, $1\le j\le
k\le g$, where $\omega\omega_i:=\omega_{\one_i}\omega_{\two_i}$, $i\in
I_M$. Set $d\tau_i:=d\tau_{\one_i\two_i}$, $i\in I_M$.

\vskip 6pt

\begin{cor}\label{threla}Fix $(p_1,\ldots,p_g)\in C^g\setminus A$, with $K(p_3,\ldots,p_g)\neq 0$.
The following $(g-2)(g-3)/2$ linearly independent relations
\begin{equation}\label{relai}\sum_{j=1}^M
C^\omega_{ij}\omega\omega_j(z):={\kuno[\sigma]\over\kdue[v]}^{g+1}F(p,x){V_{\one_i\two_i}(p_3,\ldots,p_g,x_1,\ldots,x_{2g-1},z)\over
\deltadiv\bigl(\sum_{1}^{2g-1}x_j+\b\bigr)}=0 \ ,\end{equation} $N+1\le i\le M$,
where
$$
F(p,x):=c'_{g,2}{\prod_{_{j<k}^{j,k=3}}^gE(p_j,p_k)^{g-4}\prod_{_{j\neq
\one_i}^{j=3}}^gE(p_{\one_i},p_j) \prod_{_{j\neq
\two_i}^{j=3}}^gE(p_{\two_i},p_j)\over
\prod_{j=1}^{2g-1}(\sigma(x_j)\prod_{k=3}^gE(x_j,p_k)
\prod_{k=j+1}^{2g-1}E(x_j,x_k))}\ ,
$$
with $c'_{g,2}$ is the integer constant defined in Lemma
\ref{thcombvi}, hold for all $z\in C$. Furthermore, $C^\omega_{ij}$ are independent
of $p_1,p_2,x_1,\ldots,x_{2g-1}\in C$ and correspond to the
coefficients defined in \eqref{leC} (with $\eta_i$ replaced by
$\omega_i$, $i\in I_g$) or, equivalently, in \eqref{leCi}.
Applying the Kodaira-Spencer map to Eq.\eqref{relai} yields the
$(g-2)(g-3)/2$ linear relations, $N+1\leq i\le M$,
\begin{equation}\label{finale}\sum_{j=1}^MC^\omega_{ij} d\tau_j=0\ .\end{equation}\end{cor}

\vskip 6pt

\noindent {\sl Proof.}  Consider the identities \begin{equation}\label{ratttio}{\det\nolimits_I
\sigma\sigma(x_1,\ldots,x_{2g-1},z) \over \det
v(p_3,\ldots,p_g,x_1,\ldots,x_{2g-1})}=0\ ,\end{equation} $I:=I_{M,2}\cup \{i\}$,
$N+1\le i\le M$. Upon applying Lemma \ref{thcombvi}, with $n=2$, and
Eq.\eqref{dettheta} to the numerator and Eq.\eqref{detthetaii} to the
denominator of \eqref{ratttio}, Eq.\eqref{relai} follows by a trivial
computation. On the other hand, for arbitrary points
$z,y_1,\ldots,y_{g-1}\in C$,
$$S(y_1,\ldots,y_{g-1},z)={
\sum_{i=1}^g\theta_{\Delta,i}(y_1+\ldots+y_{g-1})\omega_i(z)\over
\sigma(z)\prod_{1}^{g-1}E(z,y_i)}\ .$$ Replacing each $S$ in
$V_{\one_i\two_i}(p_3,\ldots,p_g,x_1,\ldots,x_{g-1},z)$ that depends on $z$ by the previous expression, the dependence on $z$ only enters through
$\omega_i\omega_j(z)$ and the relations \eqref{relai} can be expressed in
the form of \eqref{corol}. Eq.\eqref{finale} follows by \eqref{relai} and \eqref{kodairaspecncer}. \hfill$\square$

\vskip 6pt

\begin{theorem}\label{thdexi}If $p_1,\ldots,p_g\in C$ are $g$ pairwise distinct points
such that $K(p_3,\ldots,p_g)\neq 0$, then
$$ \Xi_i:=
\sum_{j=1}^MX^\omega_{ji}d\tau_j\ ,
$$
$i\in I_N$, with (see
Eq.\eqref{leX})
$$X^\omega_{ij}={\theta_{\Delta,\one_j}(\a_{\one_i})
\theta_{\Delta,\two_j}(\a_{\two_i})+\theta_{\Delta,\one_j}(\a_{\two_i})
\theta_{\Delta,\two_j}(\a_{\one_i})\over
(1+\delta_{\one_j\two_j})\sum_{l,m}\theta_{\Delta,l}(\a_{\one_i})
\theta_{\Delta,m}(\a_{\two_i})\omega_l(p_{\one_i})\omega_m(p_{\two_i})}\
,$$ $i,j\in I_M$, is a $\Gamma_g$-invariant basis of
$T^*_C\hat\T_g$.\end{theorem}

\vskip 6pt

\noindent {\sl Proof.}
It follows by \eqref{vXww} and  \eqref{kodairaspecncer} that
$$k(v_j)={1\over2\pi i}\sum_{k=1}^M X_{kj}^\omega d\tau_k \ ,
$$ $j\in I_N$, so that the differentials $\Xi_j:=2\pi i\,k(v_j)$,
 $j\in I_N$, are linearly independent. Furthermore, since $\{v_i\}_{i\in I_N}$ is independent of the
choice of a symplectic basis of $H_1(C,\ZZ)$, such differentials are invariant
under \eqref{modull}. \hfill$\square$

\vskip 6pt

\begin{theorem}\label{brillnoethermatrix}Let $\delta:=\left[^{\delta'}_{\delta''}\right]\in \ZZ_2^{2g}$ be a non-singular even theta characteristic
and $q\in C$ a point such that
$p_1+\ldots+p_g-q-\Delta=\delta''+\tau\delta' \mod
\ZZ^g+\tau\ZZ^g$, for some pairwise distinct points
$p_1,\ldots,p_g\in C$. Then $\forall x_1,\ldots,x_g\in C$ the determinant of the Brill-Noether matrix $\omega_i(x_j)$ is given by
\begin{align}\det \omega_i(x_j) & ={\theta[\delta](0)^g\prod_1^g
E(q,p_i)\over\det \theta_i[\delta](p_j-q)}\prod_1^g S_{\delta}(x_i,q)
\notag \\
& \cdot  {\deltadiv(\sum_1^g x_i-q)\prod_{i<j}
E(x_i,x_j)E(p_i,p_j)\over\theta[\delta](0)\prod_{i,j} E(x_i,p_j)}\ .
\end{align}\end{theorem}

\vskip 6pt

\noindent {\sl Proof.} By the Fay identity \eqref{Fay}
$$\det \lambda_i(x_j)=\prod_1^g S_{\delta}(x_i,q)
{\theta_\Delta(\sum_1^g x_i-q)\prod_{i<j} E(x_i,x_j)
E(p_i,p_j)\over\theta[\delta](0)\prod_{i,j} E(x_i,p_j)}\ ,$$ so that
by \eqref{dettheta}
$$\kuno[\lambda]={\sigma(q)\prod_1^g\theta[\delta](x_i-q)\prod_{i<j} E(p_i,p_j)\over
\theta[\delta](0)^{g+1}\prod_1^g\sigma(x_i)\prod_{i,j} E(x_i,p_j)} \
. $$
The theorem then follows by
$$\kuno[\omega]={\theta[\delta](0)^g\prod_1^g E(q,p_i)\over\det
\theta_i[\delta](p_j-q)}\kuno[\lambda] \ .$$ \hfill$\square$

\vskip 6pt

\begin{remark}\label{import} Theorem \ref{brillnoethermatrix} answers to the question, that goes back to Klein, of expressing the determinant of the
Brill-Noether matrix $\omega_i(p_j)$ in terms of theta functions and prime forms. The previously known expression is just
Eq.\eqref{dettheta}, with $\phi^1_i$ replaced by $\omega_i$, $i\in I_g$, as
given by Fay in Corollary 2.17 of \cite{jfayy}
and, more recently, in
Corollary 1.4 of \cite{FayMAM}.
On the other hand, note that such an
expression needs the section  \cite{jfayy}
$\sigma(z)=\exp-\sum_{i=1}^g\oint_{\alpha_i}\omega_i(w)\log E(w,z)$,
$z\in C$. Due to the term $\log E(w,z)$ the integral is not defined
on $C$, a problem that reflects in the constant $\kuno[\omega]$. In
\cite{FayMAM} Fay provided a definition of $\sigma$ on a fundamental
domain on the upper half-plane.
In spite of that, it is immediate to check that
$$
{\kuno[\omega]^{(2n-1)^2}\over \kenne[\phi^n] }\ ,
$$
is a well defined constant corresponding to the Mumford form of weight $n$ \cite{FayMAM}.
As will be illustrated in a forthcoming paper, the Mumford forms define vector-valued Teichm\"uller modular forms. In particular,
within the framework of the present work, there are basic relationships among such ratios, the discriminants of canonical curves
and modular forms, that in the case $g=4$ involve the Schottky-Igusa
form and the product of Thetanullwerte. This seems to be the lowest
genus case of a general relationship between Mumford forms and the
higher genus analogue of the Schottky-Igusa form, i.e. the modular
forms vanishing on the locus of Jacobians, a crucial issue in the
Schottky problem.\end{remark}

\subsection{Proofs of the main theorems}

\noindent{\it Proof of Theorem \ref{main}}

\vskip 6pt

\noindent Fix $i,j$, $3\le i<j\le g$, and choose $p_1,p_2$ in such a
way that $\{\sigma_i\}_{i\in I_g}$ is a basis of $H^0(K_C)$.
Applying Lemma \ref{thcombvi}, with $n=2$, to the identities of Theorem
\ref{thlemma} yields to the identities of  Theorem \ref{main} upon replacing
$\{\sigma_i\}_{i\in I_g}$ by $\{\omega_i\}_{i\in
I_g}$. \hfill$\square$

\vskip 6pt

To prove Theorem \ref{szeghimatr} we first need a preparatory lemma
introducing a (not normalized) Petri's like basis
$\lambda_1,\ldots,\lambda_g$ of $H^0(K_C)$, expressed in terms of
Szeg\"o kernels (defined in Eq.\eqref{szegodef}).

\begin{lem}\label{szzegibasis}Let $\delta:=\left[^{\delta'}_{\delta''}\right]\in \ZZ_2^{2g}$ be a non-singular even theta characteristic
and $q\in C$ a point such that
$A(p_1+\ldots+p_g-q)-\Delta=\delta''+\tau\delta' \mod
\ZZ^g+\tau\ZZ^g$, for some pairwise distinct points
$p_1,\ldots,p_g\in C$. The set $\{\lambda_i\}_{i\in I_g}$, where
\begin{equation}\label{szeghidef}\lambda_i(z):= S_\delta(z,q)S_\delta(z,p_i)\ , \end{equation} $i\in
I_g$, with $S_\delta(\cdot,\cdot)$ the Szeg\"o kernel, is a basis of
$H^0(K_C)$ such that
$$\lambda_i(p_j)=0\quad \hbox{\rm for}\quad 1\le i\neq j\le g\ .
$$
Furthermore, if $\phi$ is an arbitrary holomorphic section of
$L_\delta^n$, with $n>2$ an arbitrary integer, with divisor
$(\phi)=p_1+\ldots+p_{N_{n/2}}$, $N_{n/2}=n(g-1)$, for some pairwise
distinct points $p_1,\ldots,p_{N_{n/2}}\in C$, then
$$\varphi_i(z):=\phi(z) S_\delta(z,p_i)\ , $$ $i\in
I_{N_{n/2}}$, span $H^0(L_\delta^{n+1})$, and
$$\varphi_i(p_j)=0\quad\hbox{\rm for}\quad 1\le i\neq j\le N_{n/2}\ .
$$\end{lem}

\noindent{\sl Proof.} For any $i\in I_g$, the possible poles of
$\lambda_i$ at $p_i$ and $q$ cancel by
$S_\delta(p_i,q)=S_\delta(q,p_i)=0$. The only non-trivial statement
is that $\varphi_1,\ldots,\varphi_{N_n}$ are linearly independent
(the proof of linear independence of $\lambda_1,\ldots,\lambda_g$ is
analogous). Assume that $\sum_{1}^{N_n}c_i\varphi_i=0$ for some
$c_1,\ldots,c_{N_n}\in\CC$. By evaluating such an identity at
$p_1,\ldots,p_{N_n}$, it follows that $c_i\varphi_i(p_i)=0$ for all
$i\in I_{N_n}$. To prove linear independence, it is sufficient to
show that $\varphi_i(p_i)\neq 0$ for all $i\in I_{N_n}$. Suppose
this is false; then, $\varphi_i$ vanishes at $p_1,\ldots,p_{N_n}$
and $\varphi_i/\phi$ is a holomorphic section of $L_\delta$, which
is absurd since $\delta$ is non-singular.\hfill$\square$

\vskip 6pt

\noindent{\it Proof of Theorem \ref{szeghimatr}}

\vskip 6pt

\noindent Let us remind why the points $p_1,\ldots,p_g$ are uniquely determined (up to the ordering) by the condition
$p_1+\ldots+p_g-q-\Delta= \delta''+\tau\delta'$. In fact, if a different set of points $q_1,\ldots,q_g$ satisfies
the same condition, then $A(p_1+\ldots+p_g)\in W^1_g$, so that $\theta[\delta](0)=\deltadiv(p_1+\ldots+p_g-q)=0$ and $\delta$ is a singular spin structure,
contrary to the assumptions.
Define $\lambda_1,\ldots,\lambda_g\in H^0(K_C)$ as in Eq.\eqref{szeghidef}.
For any $1\le i<j\le g$,
$\lambda_i\lambda_j$ is a
holomorphic section of $K_C^2\otimes\O(-p_1\-\ldots-p_g)$. Since
$h^0(K_C^2\otimes\O(-p_1\-\ldots-p_g))=2g-3$, this implies that the $2g\times g(g-1)/2$-dimensional matrix with columns
\begin{equation}\label{szeghicols}\left(\begin{matrix}\lambda_i\lambda_{j}(x_1) \\
\vdots \\ \lambda_i\lambda_{j}(x_{2g})\end{matrix}\right)_{1\le i<j\le g}\ ,\end{equation}  for arbitrary
$x_1,\ldots,x_{2g}\in C$, has rank at most $2g-3$. Note that if $p_1,\ldots,p_g$ are pairwise distinct and $K(p_3,\ldots,p_g)\neq 0$, then
$\lambda\lambda_{g+1},\ldots,\lambda\lambda_N$ span
$H^0(K_C^2\otimes\O(-p_1\-\ldots-p_g))$ so that, in this case, the
$2g-3$ columns with $i\in\{1,2\}$ are linearly independent. Let us prove that
such a matrix has the same rank as the matrix of Theorem
\ref{szeghimatr}. We can assume that $p_1,\ldots,p_g$ are pairwise distinct; the general case is proved by applying the same argument to a maximal subset of pairwise distinct points in $\{p_1.\ldots,p_g\}$. By replacing each $\lambda_1,\ldots,\lambda_g$ in \eqref{szeghicols} by their expression in terms of Szeg\"o kernels and
neglecting the factor $S_\delta(x_i,q)^2$ appearing in each entry of
the $i$-th row, $i\in I_{2g}$, the rank of the matrix \eqref{szeghicols} is
the same as for the matrix
$$\left(\begin{matrix}S_\delta(p_i,x_1)S_\delta(p_j,x_1) \\ \vdots \\
S_\delta(p_i,x_{2g})S_\delta(p_j,x_{2g}) \end{matrix}\right)_{1\le i<j\le g}\ .
$$  Each entry in such a matrix is of the form
$S_\delta(a,x)S_\delta(b,x)$, so that it can be replaced by the RHS
of the following formula (an immediate generalization of Eq.(38) of
\cite{jfayy}),
$$S_\delta(a,x)S_\delta(b,x)= S_\delta(b,a)\Bigl(\omega_{a-b}(x)+\sum_{i=1}^g
\omega_i(x)\partial_i\log\theta[\delta](a-b)\Bigr)\ ,
$$
which holds for arbitrary points $a,b,x\in C$. For each $1\le i<j\le g$,
$S_\delta(p_i,p_j)$, $i=g+1,\ldots,N$ multiplies all the
entries of the corresponding column, so that it can be dropped. Finally,
integrate $x_i$ and $x_{i+g}$, for all $i\in I_g$, over $\alpha_i$
and, respectively, $\beta_i$. By \eqref{ominta} and \eqref{omintb}, one obtains
the matrix of Eq.\eqref{szeghione}.

\noindent To prove Eq.\eqref{szeghitwo}, note that
$A(p_i-p_j)+A(p_j-p_k)+A(p_k-p_i)=0$
for each $1\le i,j,k\le g$, so that, by taking suitable
linear combinations of the columns of the matrix of Eq.\eqref{szeghione},
we obtain that the matrix
$$\left(\begin{matrix} \vec\nabla\log \theta[\delta](p_{1}-p_{2})) & \ldots & \vec\nabla\log \theta[\delta](p_{1}-p_{g}))&
(\vec f(p_1,p_i,p_j))_{1<i<j\le g}
\\ A(p_1-p_2) & \ldots & A(p_1-p_g) & 0\end{matrix}\right)\ ,$$ has rank less than $2g-2$. For $p_1,\ldots,p_g$ in general position, the first $2g-3$
columns of the matrix in Eq.\eqref{szeghione} are linearly independent and
this implies Eq.\eqref{szeghitwo}.

Finally, consider the $g\times (g-1)$ matrix $(a_{ij})$ obtained as
a product of the matrix of Eq.\eqref{szeghithree} by the matrix
$\omega_i(x_j)$, with $x_1,\ldots,x_g$ arbitrary points. For
$i=1,\ldots,g-2$, we have $a_{ij}=\varphi_i(x_j)$ where, by
construction, $\varphi_i$ is an element of $H^0(K_C\otimes\O(-2q))$.
Since $h^0(K_C\otimes\O(-2q))=g-2$ and since the first $g-2$ columns
of the matrix of Eq.\eqref{szeghithree} are linearly independent, it
follows that $\varphi_1,\ldots,\varphi_{g-2}$ span
$H^0(K_C\otimes\O(-2q))$. On the other hand, the last column of
$(a_{ij})$ is given by
$a_{j\,g-1}=\sum_i\theta_i[\nu](0)\omega_i(x_j)=h_\nu(x_j)^2$,
$j\in I_g$, where $h_\nu$ is the unique section (up to
normalization) of $H^0(L_\nu)$. Here, $L_\nu$ is the line bundle of
degree $g-1$ with $L_\nu^2=K_C$ and associated to the
theta characteristic $\nu$, so that $(h_\nu)=\d_{g-1}$. Since $q$
is in the support of $\d_{g-1}$, $h_\nu^2\in H^0(K_C\otimes\O(-2q))$
and therefore it must be a linear combination of
$\varphi_1,\ldots\varphi_g$. In particular, ${\rm rk}(a_{ij})<g-1$
for all $x_1,\ldots,x_g$ and Eq.\eqref{szeghithree} follows by noting that
$\omega_i(x_j)$ is non-singular for general choices of
$x_1,\ldots,x_g$. \hfill$\square$


\begin{thebibliography}{99}

\bibitem{ottimo}E.~Arbarello, M.~Cornalba, P.~A.~Griffiths and
J.~Harris, {\it Geometry of Algebraic Curves}, Vol.1, Springer,
Berlin, Heidelberg, New-York, Tokyo, (1985).

\bibitem{ArbDeConc}E.~Arbarello and C.~De~Concini, On a set of
equations characterizing Riemann matrices, {\it Ann. Math.} {\bf
120} (1984), 119-140.

\bibitem{Arbar}E.~Arbarello and C.~De Concini, Another proof of a
conjecture of S.P. Novikov on periods of abelian integrals on
Riemann surfaces, {\it Duke\ Math.\ Journal\ } {\bf 54} (1987),
163-178.

\bibitem{ArbHarr}E.~Arbarello
and J.~Harris, Canonical curves and quadrics of rank $4$, {\it
Compos. Math.} {\bf 43} (1981),
145-179.

\bibitem{ArKrMar}E.~Arbarello, I.~M.~Krichever and G.~Marini,
Characterizing Jacobians via flexes of the Kummer variety, {\it
Math. Res. Lett.} {\bf 13} (2006), 109-123.

\bibitem{ArbSern} E.~Arbarello and E.~Sernesi, Petri's approach
to the study of the ideal associated to a special divisor, {\it
Invent. Math. } {\bf 49} (1978), 99-119.

\bibitem{BL}C.~Birkenhake and H.~Lange, {\it Complex Abelian Varieties}, Second edition. Grundlehren der Mathematischen Wissenschaften,
302. Springer-Verlag, Berlin, 2004.

\bibitem{Sebastian}S.~Casalaina-Martin,
Singularities of the Prym theta divisor, {\it Ann. \ of \ Math.\ } {\bf 170} (2009), 162-204.

\bibitem{DHokerQP} E.~D'Hoker and
D.~H.~Phong, Two-loop superstrings. IV: The cosmological constant
and modular forms, {\it Nucl.\ Phys.\ } B {\bf 639} (2002), 129-181.

\bibitem{donagifirst}R.~Donagi, Non-Jacobians in
the Schottky loci, {\it Ann. Math.} {\bf 126}  (1987),
193-217.

\bibitem{donagisecond}R.~Donagi, Big Schottky, {\it Invent.
Math.} {\bf 89}  (1987),  569-599.

\bibitem{donagi}R.~Donagi, The
Schottky problem, {\it Theory of moduli} (Montecatini Terme, 1985),
84-137, Springer Lecture Notes in Math. {\bf 1337},
1988.

\bibitem{dubrovin}B.~Dubrovin, Theta
functions and non-linear equations, {\it Russ.\ Math.\ Surv.\ } {\bf
36} (1981), 11-92.

\bibitem{Farkas}H.~Farkas,
Vanishing theta nulls and Jacobians, {\it The geometry of Riemann
surfaces and abelian varieties}, 37-53, Contemp. Math. {\bf 397},
2006.

\bibitem{jfayy}J.~Fay, {\it Theta Functions on Riemann
surfaces}, Springer Lecture Notes in Math. {\bf 352}, 1973.

\bibitem{FayMAM}J.~Fay, Kernel functions, analytic torsion and
moduli spaces, {\it Mem.\ Am.\ Math.\ Soc.\ } {\bf 96} (1992).

\bibitem{gunning}R.~C.~Gunning, Some curves in
abelian varieties, {\it Invent.\ Math.\ } {\bf 66} (1982), 377-389.

\bibitem{gunningtwo}R.~C.~Gunning, Some
identities for abelian integrals, {\it Amer.\ J.\ Math.} {\bf 108}
(1986), 39-74.

\bibitem{grushevsky}S.~Grushevsky, The Schottky problem, Proceedings of ``Classical Algebraic Geometry Today" workshop, MSRI 2009.

\bibitem{grushsalv}S.~Grushevsky
and R.~Salvati Manni, Jacobians with a vanishing theta-null in genus 4, {\it Israel \ J. \ Math. \ } {\bf 164} (2008) 303-315.

\bibitem{krichever}I.~M.~Krichever, Integration of non-linear equations
by methods of algebraic geometry, {\it Funct. Anal.\ Appl.\ } {\bf
11} (1) (1977), 12-26.

\bibitem{KricheverET}
  I.~M.~Krichever,
  Integrable linear equations and the Riemann-Schottky problem,
  {\it Algebraic geometry and number theory}, 497-514, Progr. Math.
  {\bf 253}, Birkh\"auser Boston, Boston, 2006.

\bibitem{Kricheverc}I.~M.~Krichever,
 Characterizing Jacobians via trisecants of the Kummer variety, {\it Ann.  \ of \ Math. \ } {\bf 172} (2010), 485-516.

\bibitem{Marini}G.~Marini, A geometrical proof of Shiota's
Theorem on a conjecture of S.P. Novikov, {\it Compos.\ Math.\ } {\bf
111} (1998), 305-322.

\bibitem{mulase}M.~Mulase, Cohomological
structure in soliton equations and Jacobian varieties, {\it J.\
Diff.\ Geom.\ } {\bf 19} (1984), 403-430.

\bibitem{mumfordd}
D. Mumford, {\it The Red Book of Varieties and
Schemes}, Springer Lecture Notes in Math. {\bf 1358}, 1999.

\bibitem{petriuno}
K.~Petri, $\rm\ddot{U}$ber die invariante darstellung
algebraischer funktionen einer ver$\rm\ddot{a}$nderlichen, {\it
Math.\ Ann.\ } {\bf 88} (1922), 242-289.

\bibitem{polishchuk}A.~Polishchuk, {\it Abelian
varieties, theta functions and the Fourier transform}, Cambridge
Univ. Press, 2003.

\bibitem{StDon}B.~Saint-Donat,
On Petri's analysis of the linear system of quadrics through a
canonical curve, {\it Math. Ann.}, {\bf 206} (1973),
157-175.

\bibitem{Schreyer}F.~O.~Schreyer, A standard basis approach to
syzygies of canonical curves, {\it J. Reine Angew. Math.}, {\bf 421}
(1991), 83-123.

\bibitem{shiota}T.~Shiota, Characterization of Jacobian varieties
in terms of soliton equations, {\it Invent.\ Math.\ } {\bf 83}
(1986), 333-382.

\bibitem{vangeemen}B.~van Geemen, Siegel modular forms
vanishing on the moduli space of curves, {\it Invent. Math.} {\bf
78} (1984) 329-349.

\bibitem{vansquare}B.~van~Geemen and G.~van~der~Geer, Kummer
varieties and the moduli spaces of abelian varieties, {\it Amer. J.
Math.} {\bf 108} (1986), 615-641.

\bibitem{welters}G.~E.~Welters, A characterization of
non-hyperelliptic Jacobi varieties, {\it Invent.\ Math.\ } {\bf 74}
(1983), 437-440.

\bibitem{Welterfg}G.~E.~Welters, A criterion for Jacobi
varieties, {\it Ann.\ Math.\ } {\bf 120} (1984), 497-504.

\end{thebibliography}
\end{document}